\documentclass{amsart}
\usepackage{graphicx}
\usepackage{amssymb}
\usepackage{epstopdf}
\usepackage{nicefrac}
\usepackage{enumerate}
\usepackage{float}
\newtheorem{thm}{THEOREM}[section]
\newtheorem{conj}[thm]{CONJECTURE}
\newtheorem{cor}[thm]{COROLLARY}
\newtheorem{defn}[thm]{DEFINITION}

\newtheorem{lemma}[thm]{LEMMA}
\newtheorem{prob}[thm]{PROBLEM}
\newtheorem{prop}[thm]{PROPOSITION}
\newtheorem{quest}[thm]{QUESTION}

\newcommand{\ds}{\displaystyle}
\newcommand{\sop}{{\bf Proof:} } % replaced by command \proof
\newcommand{\eop}{\;\;\; \Box} % replaced by command \endproof - still used inside theorem environment

\newcommand{\F}{{\mathcal F}} % our favorite foliation
\newcommand{\G}{\Gamma}

\newcommand{\e}{{\epsilon}} % some choice
\newcommand{\SOq}{{\mathbf{SO}(q)}}
\newcommand{\SO}{{\mathbf{SO}}}
\newcommand{\SL}{{\mathbf{SL}}}
\newcommand{\whF}{{\widehat{\mathcal F}}}

\newcommand{\eb}{{\overline{\exp}}}
\newcommand{\pb}{{\overline{\phi}}}

\newcommand{\ve}{\varepsilon}

\newcommand{\tn}{{\, \|\hspace{-.74pt}|\,}}

\newcommand{\wty}{{\widetilde{y}}}

\newcommand{\wtL}{{\widetilde L}}

\newcommand{\wttau}{{\widetilde{\tau}}}

\newcommand{\wtGamma}{{\widetilde{\Gamma}}}

\newcommand{\wtcP}{{\widetilde{\mathcal P}}}

\newcommand{\whg}{{\widehat{g}}}
\newcommand{\whh}{{\widehat{h}}}
\newcommand{\whm}{{\widehat{m}}}
\newcommand{\whx}{{\widehat{x}}}
\newcommand{\whz}{{\widehat{z}}}

\newcommand{\whmu}{{\widehat \mu}}
\newcommand{\whalpha}{{\widehat \alpha}}

\newcommand{\whC}{{\widehat C}}

\newcommand{\whL}{{\widehat L}}
\newcommand{\whM}{{\widehat M}}
\newcommand{\whN}{{\widehat N}}
\newcommand{\whS}{{\widehat S}}

\newcommand{\whmK}{{\widehat \mK}}

\newcommand{\bB}{{\bf B}}

\newcommand{\mB}{{\mathbb B}}
\newcommand{\mC}{{\mathbb C}}
\newcommand{\mD}{{\mathbb D}}
\newcommand{\mE}{{\mathbb E}}
\newcommand{\mF}{{\mathbb F}}
\newcommand{\mG}{{\mathbb G}}
\newcommand{\mK}{{\mathbb K}}

\newcommand{\mQ}{{\mathbb Q}}
\newcommand{\mR}{{\mathbb R}}
\newcommand{\mS}{{\mathbb S}}
\newcommand{\mT}{{\mathbb T}}
\newcommand{\mV}{{\mathbb V}}
\newcommand{\mZ}{{\mathbb Z}}

\newcommand{\cA}{{\mathcal A}}
\newcommand{\cB}{{\mathcal B}}

\newcommand{\cH}{{\mathcal H}}

\newcommand{\cL}{{\mathcal L}}

\newcommand{\cO}{{\mathcal O}}
\newcommand{\cP}{{\mathcal P}}

\newcommand{\cS}{{\mathcal S}}

\newcommand{\fG}{\mathfrak{G}}

\newcommand{\fM}{{\mathfrak{M}}}

\newcommand{\fT}{{\mathfrak{T}}}

\newcommand{\vp}{{\varphi}}

\begin{document}

\title{Embedding solenoids in foliations}

\thanks{2000 {\it Mathematics Subject Classification}. Primary 57R30, 37C55, 37B45; Secondary 53C12 }

\author{Alex Clark}
\thanks{AC supported in part by  EPSRC grant EP/G006377/1}
\address{Alex Clark, Department of Mathematics, University of Leicester, University Road, Leicester LE1 7RH, United Kingdom}
\email{adc20@le.ac.uk}

\author{Steven Hurder}
\thanks{SH supported in part by  NSF Grant 0406254}

\address{Steven Hurder, Department of Mathematics, University of Illinois at Chicago, 322 SEO (m/c 249), 851 S. Morgan Street, Chicago, IL 60607-7045}
\email{hurder@uic.edu}
\thanks{Final version date: April 27, 2011}

\date{}

% \subjclass{Primary 57R30, 37C55, 37B45; Secondary 53C12}

\keywords{}

\maketitle

\section{Introduction} \label{sec-intro}

A  \emph{ foliated space $\fM$ of dimension $n$} is a continuum
which has a local product structure:
  every point   $x \in \fM$ has an open neighborhood $U_x \subset \fM$ homeomorphic
to an open subset   $U_x \subset \mR^n$ times a compact metric space $\fT_x$. The space $\fT_x$   is called the \emph{local transverse model}, and the pre-images of the slices $U_x \times \{w\}$ for $w \in \fT_x$ are called   \emph{plaques}.  The leaves of the foliation $\F$ of $\fM$ are the
maximal connected components with respect to the fine topology on
$\fM$ induced by the plaques of the local product structure.
This concept has deep roots in continua theory, with modern treatments in the spirit of foliated manifolds given in \cite{CandelConlon2000,MS2006}. A foliated space is  \emph{transitive} if there exists at least one leaf of $\F$ which is dense in $\fM$,  and \emph{minimal} if every leaf of $\F$ is dense in $\fM$.

A \emph{matchbox manifold} is a foliated space $\fM$ such that all of the
local transverse models $\fT_x$ are $0$-dimensional, or equivalently are totally disconnected \cite{AO1991,AO1995}.
 In this case, the leaves of $\F$ are simply the path components of $\fM$. Intuitively, a
$1$-dimensional matchbox manifold $\fM$ has local coordinate charts
$U$ which are homeomorphic to a ``box of matches.'' This intuitive description is extended to higher dimensional leaves with a stretch of the imagination.

One source of examples of transitive foliated spaces is to consider a   compact manifold $M$ and leaf $L$ of a foliation $\F_M$ of $M$, then let $\fM = \overline{L}$ denote the topological closure of $L$ in $M$. The topological space $\fM$ inherits the structure of a foliated space, where $\F = \F_M \mid \fM$. Zorn's Lemma implies that each such closure  $\fM = \overline{L}$ also contains at least one minimal foliated space.
The question we consider in this note is   the  converse of these examples:
\begin{prob}\label{prob-embed}
Let $\fM$ be a minimal matchbox manifold of dimension $n$. When does there exists a $C^r$-foliation $\F_M$ of a compact manifold $M$ and a foliated topological embedding  $\iota \colon \fM \to M$ realizing $\fM$ as a   foliated subset?
\end{prob}

 The embedding problem is fundamentally different for the cases of differentiability classes  $r=0$ and  $r \geq 1$.
 In the case of $C^0$-foliations, Clark and Fokkink \cite{ClarkFokkink2004} obtain some general embedding results for codimension $q \geq 2$.   There are no known obstructions to the existence of some   embedding of a minimal matchbox manifold $\cP$ into some $C^0$-foliated compact manifold.

 What makes the embedding problem particularly interesting for $r \geq 2$, is that an embedding of $\fM$ as  a minimal set of a $C^r$-foliation $\F_M$ of a compact manifold $M$ is  equivalent to a delicate question concerning the possible holonomy pseudogroups of $C^r$-foliations \cite{Haefliger2002}. Solutions to this problem have applications to the study of the homotopy type of the Haefliger classifying space $B\G^r_n$ for $C^r$-foliations \cite{Haefliger1984, Hurder2009,Hurder2011b}.

One source of examples of matchbox manifolds is provided by the
space of tilings associated to a given quasi-periodic tiling of
$\mR^n$, as in \cite{AP1998,FHK2002,Ghys1999}, and also by the more
general case of $G$-solenoids as introduced by Benedetti and
Gambaudo in \cite{BG2003}. For a few classes of quasi-periodic
tilings of $\mR^n$, the codimension one canonical cut and project
tiling spaces ~\cite{FHK2002}, it is known that the associated
matchbox manifold is a minimal set for a $C^1$-foliation of a torus
$\mT^{n+1}$, where the foliation is a generalized Denjoy example.

Another class of examples of naturally occurring matchbox manifolds
are provided by the ``Williams  solenoids'', as introduced in
\cite{Williams1967,Williams1974} to describe the attractors of
certain Axiom A attractors. While Williams finds embeddings of these
matchbox manifolds as attractors of diffeomorphisms of smooth
manifolds such that the path components of these matchbox manifolds are
leaves of the unstable foliation of the corresponding
diffeomorphism, it is unknown which of the Williams solenoids can be
embedded as minimal sets for foliations of the entire ambient manifold.

The embedding problem has been solved for the case of  Markov
minimal sets \cite{CantwellConlon1988, Matsumoto1988} associated to
pseudogroups of local $C^r$-diffeomorphisms of the line $\mR$.
Cantwell and Conlon \cite{CantwellConlon2002} gave an explicit
construction for realizing any such matchbox manifold  as a minimal
set for a codimension-one $C^r$ foliation of a compact manifold.

An important   special case of minimal matchbox manifolds   are the solenoids studied by McCord in \cite{McCord1965}, spaces which admit a presentation as an inverse limit of closed manifolds. This concept greatly generalizes the classical  solenoids  introduced by Vietoris  \cite{Vietoris1927},   which are modeled on the circle $\mS^1$ .
\begin{defn}  \label{def-solenoid}
Assume there is given a collection $\cP = \{ p_{\ell} \colon M_{\ell} \to M_{\ell -1} \mid \ell \geq 0\}$, where each $M_{\ell}$ is a connected compact manifold without boundary of dimension $n$, and
$p_{\ell} \colon M_{\ell} \to M_{\ell -1}$  is a covering map of degree $d_{\ell} > 1$. Then the    inverse limit topological space
\begin{equation}\label{eq-invlim}
 \cS = \cS_{\cP} \equiv \lim_{\leftarrow} ~ \{ p_{\ell} \colon M_{\ell} \to M_{\ell -1}\}
\end{equation}
is said to be a \emph{solenoid} with base space  $M_0$, and the collection $\cP$ is said to be a \emph{presentation} for $\cS$.

If each   covering map
$\pi_{\ell} = p_1 \circ \cdots \circ  p_{\ell} \colon M_{\ell} \to M_0$ is a normal covering for all $\ell \geq 1$, then $\cS$ is said to be a \emph{ McCord solenoid}  \cite{McCord1965, Schori1983, ClarkFokkink2004}.
A presentation $\cP$ is \emph{oriented} if all manifolds $M_{\ell}$  are oriented, and the maps $p_{\ell}$ are orientation-preserving.
\end{defn}
A \emph{Vietoris  solenoid} is an inverse limit space  as in Definition~\ref{def-solenoid}  with  
base space  $M_0 = \mS^1$.

It is fundamental in Definition~\ref{def-solenoid}  that the bonding maps are covering maps, as this implies
that  $\cS$ admits a fibration $\pi  \colon \cS \to M_0$ whose fibers are totally disconnected.
(In the case of McCord solenoids, the fibers are Cantor groups.) This implies that a solenoid $\cS$
is a matchbox manifold, where the leaves of  $\F_{\cS}$  are the path components of $\cS$. Each leaf of  $\F_{\cS}$ is a  smooth manifold of  dimension $n$, and  the restriction of $\pi$ to a leaf is a smooth covering map of $M_0$. If $\cS$ is defined by an oriented presentation   $\cP$, then $\cS$ is an oriented   matchbox manifold.

  \begin{defn}  \label{def-MMS}
A matchbox manifold $\fM$ is a solenoid if there exists a collection of proper covering maps
$\cP  \equiv \{ p_{\ell} \colon M_{\ell} \to M_{\ell -1} \mid \ell \geq 0 \}$   and a homeomorphism $\vp \colon \fM \cong \cS_{\cP}$. We then say that $\{\cP, \vp\}$ (or simply $\cP$ by abuse of notation) is a \emph{presentation} for $\fM$.
\end{defn}

The purpose of this work is to study the embedding Problem~\ref{prob-embed} for solenoids.   

A matchbox manifold $\fM$ is said to be homogeneous if its group of homeomorphisms ${\bf Homeo}(\fM)$ acts transitively. McCord showed in \cite{McCord1965} that every McCord solenoid is homogeneous. The converse was proven by the authors in \cite{ClarkHurder2010}, that an orientable homogeneous matchbox manifold is a McCord solenoid, in the sense of Definition~\ref{def-MMS}. More generally, the same paper shows that if the dynamics of the foliation of an oriented  matchbox manifold $\fM$ is equicontinuous and has no germinal holonomy, then $\fM$ is a solenoid,  though not necessarily homogeneous. Thus, the solenoids are a natural dynamically-defined class of matchbox manifolds to study.

A fundamental point in our analysis of the embedding problem for solenoids, is  that a solenoid admits many presentations, all with   homeomorphic inverse limits \cite{FO2002, Mardesic2000,
McCord1965,Rogers1970, RogersTollefson1971, Schori1983}. This is
well-known for the traditional solenoids of dimension one, where
$M_0 = \mS^1$, and the homeomorphism class is determined by the list
of primes  with   multiplicity which divide  the degrees of the
covering maps $p_{\ell}$. More generally, one has the following
homeomorphism result, which is a special case of the general fact
that the inverse limit of a cofinal sequence of an inverse system is
homeomorphic to the inverse limit of the original inverse system~
\cite{CordierPorter1989, Mardesic2000}:
 \begin{prop}\label{prop-subtower}
 Let $\cP  \equiv \{ p_{\ell} \colon M_{\ell} \to M_{\ell -1} \mid \ell \geq 0 \}$    be a presentation for a solenoid $\cS$. Given a subsequence $\{0 = i_0 < i_1 < \cdots < i_{\ell} < \cdots \}$ define bonding maps
 $$p_{\ell}' \colon M_{i_{\ell}} \to M_{i_{\ell -1}} ~ , \quad  p_{\ell}' = p_{i_{\ell}} \circ   p_{i_{\ell} - 1} \circ   p_{i_{\ell} - 2}  \circ \cdots \circ p_{i_{\ell -1} + 1}$$
 Define $\cP' = \{p_{\ell}' \colon M_{i_{\ell}} \to M_{i_{\ell -1}} \mid \ell \geq 0\}$, which is said to be a
  \emph{sub-presentation} of $\cP$. Then for any sub-presentation $\cP'$ of $\cP$,    there is a natural homeomorphism $\cS_{\cP'} \cong \cS_{\cP}$.
 \end{prop}

For the  general case of an arbitrary $\fM$, the existence of a ``presentation for $\fM$'' is discussed in \cite{ClarkHurderLukina2011}, though now the notion of a presentation must allow for the model spaces $M_{\ell}$ to be branched manifolds as in the works of Williams \cite{Williams1967,Williams1974}. 

An important concept we use is the notion of  an embedding of a presentation $\cP$  of the solenoid $\cS_{\cP}$ into a foliated manifold $M$.
Given $\cP  \equiv \{ p_{\ell} \colon M_{\ell} \to M_{\ell -1} \mid \ell \geq 0 \}$, an embedding of $\cP$ 
  is a collection of homeomorphisms $\vp_{\ell} \colon M_{\ell} \to L_{\ell}$, where each $L_{\ell}$ is a leaf of $\F_M$ for $\ell \geq 0$. We also require that   there are covering maps   $\psi_{\ell} \colon L_{\ell} \to L_{\ell -1}$ which are naturally induced by the transverse geometry of $\F_M$, and an embedding   $\vp \colon \cS_{\cP} \to   M$, so that the maps $\vp_{\ell}$ induce a homeomorphism $\vp^*$ between $\cS_{\cP}$ and the inverse limit $\cS_{\F}$ of  the collection  $\cP_{\F} = \{\psi_{\ell} \colon L_{\ell} \to L_{\ell -1} \mid \ell \geq 0\}$. Thus, an embedding of the presentation $\cP$ of the solenoid $\cS_{\cP}$ is an embedding of $\cS_{\cP}$ of a very special type. 

While this type of embedding has not explicitly been studied in the literature, one finds in the construction of Markus and Meyer~\cite{MM1980} a similar concept. There, for each homeomorphism class of one-dimensional solenoids, they choose a   presentation $\cP$ of the solenoid $\cS_{\cP}$ and the authors identify in the generic Hamiltonian flow a limit set that is homeomorphic to $\cS_{\cP}$ by finding a sequence of periodic orbits of increasing period related to the presentation $\cP$ in such a way that these periodic orbits converge to a minimal set of the flow which is homeomorphic  to $\cS_{\cP}.$  However, there are no related results known to the authors for higher dimensional solenoids.

 This then leads to a refinement of the embedding Problem~\ref{prob-embed}:

 \begin{prob}\label{prob-embed2}
Let $\cS$ be an $n$-dimensional solenoid with presentation $\cP$. When does there exists a $C^r$-foliation $\F_M$ of a
compact manifold $M$ and an embedding  of $\cP$ realizing $\cS$ as the inverse limit of the associated system of leaves?
\end{prob}

Our approach to Problem~\ref{prob-embed}  in the case where $\cS$ is a solenoid, is based on the existence of a   sub-presentation   $\cP'$  with $\cS \cong \cS_{\cP'}$  where $\cP'$ is   well-chosen so that the constraints imposed by the $C^r$-differentiability requirements on an embedding of $\cP'$ are satisfied.  That is, 
given an $n$-dimensional    solenoid $\cS$ with presentation $\cP$, and so $\cS \cong \cS_{\cP}$, to solve   Problem~\ref{prob-embed} for $\cS$, it suffices by Proposition~\ref{prop-subtower}  to exhibit a sub-presentation $\cP'$ of $\cP$ and a solution to Problem~\ref{prob-embed2} for  $\cP'$.

In the case of codimension one embeddings, if $\iota \colon \cS \to
M$ is any topological embedding of a solenoid of dimension $n$  into
a compact oriented manifold of dimension $n+1$, then  Clark and
Fokkink showed in \cite{ClarkFokkink2004} that the inverse limit
(\ref{eq-invlim}) defining $\cS$ has a finiteness property for its
top degree cohomology:    there exists some $\ell_0 \gg 0$ such that
the Cech cohomology $\check{H}^n(\cS; \mR) \cong H^n(M_{\ell}; \mR)$
for all $\ell \geq \ell_0$, which is a contradiction.
Thus, in the following, one always assumes that the codimension $q \geq 2$.

   The existence of 1-dimensional solenoids as   minimal sets of smooth flows  has an extensive history in topological dynamics,  and we cite only some selected results here ~\cite{BowenFranks1976, GT1990, GST1994, Kan1986, MM1980, Smale1967, EThomas1973}. The existence is generally shown via an iterated perturbation argument, which is essentially folklore. That is, starting with a closed orbit, $M_0 \cong \mS^1$, it is modified in an open neighborhood of $M_0$ so that the perturbed flow now has a nearby closed orbit $M_1 \cong \mS^1$ which covers $M_0$ with degree $d_1 > 1$. This process is inductively repeated for all subsequent closed orbits $M_{\ell}$ with $\ell > 1$. With suitable care in the choices, the resulting flow    will  be   $C^{\infty}$ and has a minimal set $\cS$ homeomorphic to the inverse limit of the   system of closed orbits resulting from the construction. The generalization of this folklore construction to higher dimensions   requires a more precise description of the steps of the construction, as additional issues arise in higher dimensions, which are discussed in this paper.

Given a $C^r$-flow $\vp$ on a compact manifold with solenoidal minimal set $\cS_{\vp}$,  the existence of a collection of closed orbits for $\vp$   whose limit defines $\cS_{\vp}$ is one of the standard ``questions'', as such a system may or may not exist, depending on the ``geometry'' of the flow. In the formulation of the embedding Problem~\ref{prob-embed2}, such a system of closed orbits (or compact leaves)  is part of the assumption.

\section{Main Theorems} \label{sec-results}

In  this paper we give solutions to the embedding Problems~\ref{prob-embed} and \ref{prob-embed2} for a special class of McCord solenoids of dimension $n \geq 1$, those which have abelian Cantor groups as fibers.

Let $\cS$ be a McCord solenoid with presentation  $\cP = \{ p_{\ell}
\colon M_{\ell} \to M_{\ell -1} \mid \ell \geq 0\}$. Fix a basepoint
$x_0 \in M_0$. Then, inductively choose  basepoints $x_{\ell} \in
M_{\ell}$ such that $p_{\ell}(x_{\ell}) = x_{\ell -1}$ for each
$\ell \geq 1$. Set $\ds \G = \pi_1(M_0 , x_0)$, and let $\ds
\G_{\ell} = {\rm Image}\{\pi_{\ell *}  \colon \pi_1(M_{\ell},
x_{\ell}) \to \G\}$. Let $\mathfrak{K}_{\ell} = \G/\G_{\ell}$ be the
quotient set, which is a finite group by the normality assumptions.
The fiber   of   $p_* \colon \cS \to M_0$  is then identified with
the inverse limit, $\ds \mK_0  \equiv p_*^{-1}(x_0) \cong
\lim_{\leftarrow} ~ \{ \mathfrak{K}_{\ell +1} \to
\mathfrak{K}_{\ell} \}$, which is a  compact, totally-disconnected,
perfect topological group; that is, a Cantor group. Moreover, by
results of McCord \cite{McCord1965}, $\cS$ is homeomorphic to a
principle $\mK_0$-bundle over the base space $M_0$.

Now suppose that $M_0$ is a closed oriented manifold of dimension $n$ with basepoint $x_0 \in M_0$, and  that there is a surjective map $\pi \colon \G = \pi_1(M_0 , x_0) \to \mZ^k$ for some $k \geq 1$. Given a descending chain of subgroups,
 $\ds \G_{\cP} \equiv \left\{ \mZ^k = \G_0 \supset \G_1 \supset \G_2 \supset \cdots \right\}$, where each $\G_{\ell + 1} \subset \G_{\ell}$ has finite index greater than one for all $\ell \geq 0$. Set $\wtGamma_{\ell} = \pi^{-1}(\G_{\ell})$ and let $\pi_{\ell} \colon M_{\ell} \to M_0$ be the covering corresponding to   $\wtGamma_{\ell}$. Then we obtain a presentation $\wtcP$ for a McCord solenoid $\cS_{\wtcP}$ with base space $M_0$. The fiber $\mK_0$ of $p_* \colon \cS_{\wtcP} \to M_0$ is naturally isomorphic to the abelian Cantor group $\mK_0$ defined by the chain $\G_{\cP}$.

Let $\mD^q$ denote the closed unit disk in $\mR^q$. By Proposition~\ref{prop-transverse}, a $C^r$-foliation $\whF$ on  $\mT^k \times \mD^q$  determines a $C^r$-foliation $\F$ on $M \times \mD^q$ with the same dynamical properties. A compact leaf of $\whF$ lifts to a compact leaf of $\F$, and the same applies for minimal sets of $\whF$. Thus, for a solenoid $\cS_{\wtcP}$ with base $M_0$ obtained in this way, in order to construct  an embedding  of $\cS_{\wtcP}$ into a $C^r$-foliation of $M \times \mD^q$, it suffices to construct an   embedding of the   solenoid $\cS_{\cP}$ with base $\mT^k$ into a $C^r$-foliation of $\mT^k \times \mD^q$.

For any $k \geq 1$ and $r \geq 0$, and presentation $\cP$ of a McCord solenoid $\cS$ with base space $\mT^k$,  we give   explicit  constructions of  $C^r$-foliations which realize sub-presentations $\cP'$ of $\cP$. For the case of $k=1$, the construction is simply the classical method, expressed precisely in the language of flat bundles.  For dimensions $k \geq 2$, the approach using flat bundles seems fundamental, as the proofs of $C^r$-embedding for $r \geq 2$ depend upon careful analysis of the holonomy representations of the normal bundles to the compact leaves created inductively in the construction.

For the classical case of one-dimensional solenoids,  each $M_{\ell} = \mS^1$ and the bonding maps $p_{\ell} \colon \mS^1 \to \mS^1$ are orientation-preserving covering maps with degree $d_{\ell} > 1$. The foliation on the inverse limit space $\cS$   defines a minimal flow, and $\mK_0$ is homeomorphic to a   $\vec{d}$-adic Cantor set, for $\vec{d} = (d_1, d_2, \ldots)$. The solenoid $\cS$ is a compact abelian topological group, hence by Pontrjagin Duality \cite{Baer1937, Kechris2000, Pontrjagin1934}  is determined up to homeomorphism  by the list of primes (with their multiplicities) dividing the integers $d_{\ell}$. (For further discussion, see the remarks at the end of section~\ref{sec-foliations}.)  In particular, there are an uncountable number of one-dimensional solenoids which are distinct up to homeomorphism.

For the case where $k > 1$,  and a standard (affine) solenoid $\cS$  associated to a descending chain $\G_{\cP} \equiv \left\{ \mZ^k = \G_0 \supset \G_1 \supset \G_2 \supset \cdots \right\}$, the topological type of $\cS$ is determined by the chain $\G_{\cP}$. Again, there are an uncountable number of such $k$-dimensional solenoids which are distinct up to homeomorphism.
Unlike the $1$-dimensional case, the isomorphism problem for the solenoids $\cS$ resulting from such chains   is  unclassifiable in the sense of descriptive set theory \cite{Kechris2000, Thomas2001, Thomas2003, HT2006}. Thus, in contrast to the $1$-dimensional case, it is not possible to give a classification for the family of minimal sets we construct below. In all cases $k \geq 1$, our constructions yield an  uncountable number of distinct homeomorphism types of  solenoids which admit embeddings into $C^{\infty}$-foliations.

Our strongest results are for $C^0$-embedding problem. Proposition~\ref{prop-embedC0} yields the following result,  that every presentation of a solenoid with base $\mT^k$ admits an embedding into a $C^0$-foliation.
\begin{thm}\label{thm-main-top}
Let $\cP$ be a presentation of the solenoid $\cS$ over the base space $\mT^k$, and let $q \geq 2k$. Then   there exists a $C^0$-foliation $\whF$ of $\mT^k \times \mD^q$  such that:
\begin{enumerate}
\item  $\whF$  is a distal foliation, with   a smooth transverse invariant volume form;
\item $L_0 = \mT^k \times \{ \vec0\}$ is a leaf of $\whF$, and  $\whF = \F_0$ near the boundary of $M$;
\item there is an embedding of $\cP$ into the foliation $\whF$;
\item the solenoid $\cS$ embeds as  a minimal set $\whF$.
\end{enumerate}
\end{thm}

The embedding problem for solenoids into $C^1$-foliations  is the next most general case. While it is unknown whether smoothness conditions restrict the possible embeddings, new ideas must be introduced to address the smoothness conditions. In section~\ref{sec-foliations}, we introduce the ``standard representation'' for a descending chain $\G_{\cP}$ of abelian lattices, and formulate the condition (\ref{eq-C1algbounds})    on the algebraic structure of the chain $\G_{\cP}$. By careful choice of a sub-presentation $\cP'$ of the given presentation $\cP$,  Proposition~\ref{prop-embedC1} then yields:
\begin{thm}\label{thm-main-torus}
Let $\cP$ be a presentation of the solenoid $\cS$ over the base space $\mT^k$, and let $q \geq 2k$. Suppose that $\cP$ admits a sub-presentation $\cP'$ which satisfies condition  (\ref{eq-C1algbounds}).
Then   there exists a   $C^1$-foliation $\whF$ of $\mT^k \times \mD^q$   such that:
\begin{enumerate}
\item  $\whF$  is a distal foliation, with   a smooth transverse invariant volume form;
\item $L_0 = \mT^k \times \{ \vec0\}$ is a leaf of $\whF$, and  $\whF = \F_0$ near the boundary of $M$;
\item there is an embedding of $\cP'$ into the foliation $\whF$;
\item the solenoid $\cS$ embeds as  a minimal set $\whF$.
\end{enumerate}
\end{thm}

Our constructions have    applications   to the study of dynamical systems defined by foliations. In particular, we give some  ``non-stability'' results which complement the various forms of the Reeb-Thurston stability theorems for $C^1$-foliations \cite{LR1977, Stowe1983, Thurston1974}.
The generalized Reeb Stability  Theorem    \cite{Stowe1983, Thurston1974}  states that   if $L_0$ is a compact leaf of a $C^1$-foliation $\F$, such that the first cohomology group $H^1(L_0, \mV) = 0$, where $\mV$ is a vector space module over $\pi_1(L_0, x_0)$, then all leaves of $\F$ sufficiently close to $L_0$ are diffeomorphic to it.   Theorem~2 of Langevin and Rosenberg \cite{LR1977} states that if $H^1(L_0; \mR) =0$ and $L_0$ has a product open neighborhood, then if $\F'$ is a sufficiently close $C^1$-perturbation of $\F$,   there is a compact leaf $L_0'$ of $\F'$ near to $L_0$ which has a product open neighborhood.

In contrast to the assumptions for stability, Corollary~\ref{cor-uncountable} below  assumes that $H^1(L_0, \mR) \ne 0$, and shows
that there exists  a perturbation for which there exists a nearby solenoidal minimal set. Thus, the leaf $L_0$ is highly \emph{non-stable}.
This is based on our construction of $C^r$-embeddings of solenoids with base $\mT^k$ which are arbitrarily close to the product foliation:
\begin{thm}\label{thm-main-fol}
Let $L_0$ be a closed oriented manifold of dimension $n$, with $H^1(L_0, \mR) \ne 0$. Let $q \geq 2$   and $\ve> 0$.  Then there exists a $C^{\infty}$-foliation $\F$ of  $M =  L_0 \times \mD^q$  which is $\ve$-$C^{\infty}$-close to the product foliation $\F_0$, such that
 $\F$ is a volume-preserving, distal foliation, and satisfies
\begin{enumerate}
\item $L_0 = L_0 \times \{0\}$ is a leaf of $\F$;
\item $\F = \F_0$ near the boundary of $M$;
\item $\F$ has a minimal set $\cS$ which is a   solenoid with base $L_0$;
\item $\cS$ is in the closure of the compact leaves of $\F$.
 \end{enumerate}
\end{thm}

The   manifold $M =  L_0 \times \mD^q$ with foliation  $\F$ functions as a ``foliated plug'', in the spirit of the constructions of Wilson \cite{Wilson1966}, Schweitzer \cite{Schweitzer1974} and Kuperberg \cite{Kuperberg1994}. We may use it to insert  a solenoidal minimal set  into a given  foliation satisfying local conditions. This leads to the following surprising result:

\begin{cor}\label{cor-uncountable}
Let $\F_0$ be a $C^{\infty}$-foliation of codimension $q \geq 2$ on a manifold $M$. Let $L_0$ be a compact leaf with
$H^1(L_0 ; \mR) \ne 0$, and suppose that $\F_0$ is a product foliation in some saturated open neighborhood $U$ of $L_0$. Then there exists a foliation $\F_M$ on $M$ which is $C^{\infty}$-close to $\F_0$, and $\F_M$ has an uncountable set of  solenoidal minimal sets $\{\cS_{\alpha} \mid \alpha \in \cA\}$, all contained in $U$, and \emph{pairwise non-homeomorphic}. \
If $\F_0$ is a distal foliation with a smooth transverse invariant volume form, then the same holds for $\F_M$.
\end{cor}

The minimal sets of an equicontinuous foliation are submanifolds \cite{ALC2009}, hence if $\F$ has a solenoidal minimal set, then $\F$ cannot be equicontinuous. Thus,  the minimal set $\cS$ is parabolic in terms of the classification scheme of \cite{Hurder2009,Hurder2011a}.

The construction of the foliation $\F$ on $M = \mT^k \times \mD^q$
in  Theorems~\ref{thm-main-top} and \ref{thm-main-torus}   is just a special case of a
more general construction. The key technical idea is to define a type of ``plug'',
whereby an open foliated tubular neighborhood of a compact leaf is replaced by a new foliation,
which is defined using modifications of flat bundles over the leaf. In other words, the classical approach
where a periodic orbit of a flow is locally modified so that the orbits of the new flow finitely covers this orbit nearby,
is considered as the simplest case of a plug obtained from a flat bundle with finite holonomy.

Section~\ref{sec-flat} first discusses
the generalities of the construction and properties of flat vector
bundles and their modifications.  Section~\ref{sec-plug}
then shows how to use a sequence of such flat bundle modifications to obtain a sequence of
foliations $\whF_{\ell}$ on $\mT^k \times \mD^q$ which approximate
an embedding of the first $\ell$-stages of a tower of coverings of
$\mT^k$. This section is given in very geometric terms, as this
motivates the construction.  However, to prove that the sequence of
foliations $\whF_{\ell}$ converge in the uniform   $C^r$-norm,
section~\ref{sec-estimates}   reformulates the construction   in
terms of the global holonomy maps $\whh_{\ell} \colon \mZ^k \to {\rm
Diff}^r(\mD^q)$ of the foliations $\whF_{\ell}$.  This formulation
also makes the comparison of the properties of our examples  with
those of the classical examples of actions on the disk more
transparent.

The work in sections~\ref{sec-flat} to \ref{sec-estimates} assume as given  data, that there is given a sequence of   representations $\whalpha = \{ \alpha_{\ell} \colon \mZ^k \to \mQ^m \mid \ell \geq 1\}$ which determine the flat bundle perturbations at each stage.
This yields a sequence of modified foliations  $\{\whF_{\ell} \mid \ell \geq 0\}$    of $\whN_0 = \mT^k \times \mD^q$ as constructed in section~\ref{sec-plug}. The delicate issues of convergence in the $C^r$-norm in section~\ref{sec-estimates}  determine the smoothness of the limiting foliation $\whF$ of the plug $ \mT^k \times \mD^q$.

Theorems~\ref{thm-main-top} and \ref{thm-main-torus} do not specify a sequence of representations $\whalpha$; this additional data may be derived from the given data of a tower of coverings $\cP$ of the base manifold $\mT^k$. We give one example of this, the ``standard representations'' we introduce in section~\ref{sec-foliations},  as direct generalizations of the $1$-dimensional method which is presented in   section~\ref{sec-flows}.

  Section~\ref{sec-remarks} discusses   some   natural  questions concerning   solenoidal minimal sets for foliations which arise from this work.

\section{Flat bundles}\label{sec-flat}

We recall some of the basic facts about flat vector bundles, and then introduce the basic technique of using deformations of flat bundles to construct the foliated plugs which are fundamental to our constructions.

Choose a basepoint $x_0 \in L_0$, and set $\Gamma = \pi_1(L_0, x_0)$.
Let $\G$ act on the right as deck transformations of the universal cover $\widetilde{L_0} \to L_0$.

Let $\rho \colon \G \to \SOq$ be an orthogonal representation.
Then $\rho$ defines an action of $\G$   on the left as isometries of $\mR^q$, defined by  $\vec{v} \mapsto \rho(\gamma)(\vec{v})$.
Define a flat $\mR^q$-bundle by
\begin{equation}\label{eq-flat}
\mE_{\rho}^q = (\widetilde{L_0} \times \mR^q)/(\wty \cdot \gamma, \vec{v}) \sim (\wty, \rho(\gamma)( \vec{v})) ~ \longrightarrow ~ L_0
\end{equation}
For a closed path $\sigma_{\gamma} \colon [0,1] \to L_0$ with $\sigma_{\gamma}(0) = \sigma_{\gamma}(1) = x_0$
which represents the homotopy class $\gamma$, the holonomy of the bundle $\mE_{\rho}$ along the path $\sigma_{\gamma}$ is $\rho(\gamma)$, as seen via the identification in (\ref{eq-flat}).

For each $\vec{v} \in \mR^q$, let
$\wtL_{\vec{v}} = \wtL_0 \times \{\vec{v}\} $ be the leaf through $\vec{v}$ of the product foliation on $(\widetilde{L_0} \times \mR^q)$. The product foliation is $\G$-equivariant, so descends to a foliation denoted by $\F_{\rho}$, and the leaf $\wtL_{\vec{v}}$ descends to a leaf denoted by $L_{\vec{v}} \subset \mE_{\rho}^q$. The tangent bundle  $T\F_{\rho}$ defines the  integrable horizontal distribution for
the  flat bundle $\mE \to L_0$.

The action of $\G$ on $\mR^q$ via $\rho$ preserves the usual norm,
$\|\vec{x}\|^2 = x_1^2 + \cdots + x_q^2$, so for all $\e > 0$ restricts to actions on the subsets
$$
\mB^q_{\e} = \{ \vec{x} \mid\| \vec{x} \| < \e\} ~, ~
\mD^q_{\e} = \{ \vec{x} \mid\| \vec{x} \| \leq \e\} ~ , ~
\mS^{q-1}_{\e} = \{ \vec{x} \mid\| \vec{x} \| = \e\}
$$
Recall that $\mD^q = \mD^q_{1}$ and $\mB^q = \mB^q_{1}$.

Let $\mB^q_{\e, \rho} \to L_0$ (respectively, $\mD^q_{\e, \rho}$ or $\mS^{q-1}_{\e, \rho}$) denote the $\mB^q_{\e}$-subbundle (respectively, $\mD^q_{\e}$ or $\mS^{q-1}_{\e}$) of $\mE_{\rho}^q \to L_0$. The foliation $\F_{\rho}$ then restricts to a foliation on each of these subbundles.

The most familiar example is for $L_0 = \mS^1$, where $\G = \pi_1(\mS^1 , x) = \mZ$. Given any $\alpha \in \mR$, define the representation
$\rho \colon \mZ \to {\SO(2)}$ by $\rho(n) = \exp (2\pi \sqrt{-1} \, \alpha \cdot n)$.
Then $\mE^2_{\rho}$ is the flat vector bundle over $\mS^1$, such that the holonomy along the base $\mS^1$ rotates the fiber by $\rho(1) = \exp(2\pi \sqrt{-1} \, \alpha)$.

In general, the bundle $\mE_{\rho}^q \to L_0$ need not be a product vector bundle,
and may even have non-trivial Euler class when $q$ is even \cite{KT1967,Milnor1958, Sullivan1974,Wood1971}.

\begin{prop}\label{prop-trivial}
Let $\rho_t \colon \G \to \SOq$, for $0 \leq t \leq 1$, be a smooth 1-parameter family of representations such that $\rho_0$ is the trivial map, and $\rho_1 = \rho$. Then $\rho_t$ canonically defines an isomorphism of vector bundles $\mE_{\rho}^q \cong L \times \mR^q$.
\end{prop}
\sop
The family of representations defines a family of flat bundles $\mE^q_{\rho_t}$ over the product space $L \times [0,1]$.
The coordinate vector field $\partial/\partial t$ along $[0,1]$ lifts to a vector field $\widetilde{\vec{v}}$ on the product,
$\widetilde{\mE}^q = (\widetilde{L_0} \times [0,1]) \times \mR^q$. The group $\G$ acts on $\widetilde{\mE}^q$ via the action of $\rho_t$ on each slice
$\widetilde{\mE}_t^q = (\widetilde{L_0} \times \{t\}) \times \mR^q$, and the vector field $\widetilde{\vec{v}}$ is clearly $\G$ invariant, as $\rho_t$ acts on the vector space $\mR^q$ but fixes the tangent bundle to $[0,1]$. Thus, $\widetilde{\vec{v}}$ descends to a vector field $\vec{v}$ on $\widetilde{\mE}^q/\G$. The flow of $\widetilde{\vec{v}}$ on $\widetilde{\mE}^q$ preserves the fibers, so the flow of $\vec{v}$ acts via bundle isomorphisms on
$\widetilde{\mE}^q/\G \to L_0 \times [0,1]$.

The time-one flow of $\vec{v}$ defines an isotopy between the bundles $\mE^q_{\rho_0}$ and $\mE^q_{\rho_1}$, which induces a bundle isomorphism between them. The initial bundle $\mE^q_{\rho_0}$ is a product, hence the same holds for $\mE^q_{\rho_1}$.
$\eop$

The key point is that the bundle isomorphism between $\mE^q_{\rho_0}$ and $\mE^q_{\rho_1}$ is canonical, and in particular, depends smoothly on the path $\rho_t$.
In general, a representation $\rho$ need not be homotopic to the trivial representation, so the hypotheses of Proposition~\ref{prop-trivial} is very strong.

If $\rho \colon \G \to \SOq$ factors through either a free abelian group $\mZ^k$, or a free non-abelian group $\mF^k$, then one can construct an explicit isotopy $\rho_t$ of $\rho$   to the identity.  For the case where
$$ \rho = \alpha \circ \pi \colon \G \stackrel{\pi}{\longrightarrow} \mZ^k \stackrel{\alpha}{\longrightarrow} \SOq$$
we give the construction of the the isotopy $\rho_t$ in detail.
Let $m$ be the greatest integer such that $2m\leq q$. For $q = 2m$, identify $\mR^{q} \cong \mC^m$ via the map
$$\vec{x}= (x_1 , x_2, \ldots , x_{2m-1}, x_{2m}) \leftrightarrow \vec{z} = (z_1, \ldots , z_m)$$
where $z_i = x_{2i -1} + \sqrt{-1} \cdot x_{2i}$. For $q = 2m+1$, identify $\mR^q = \mR^{2m} \times \mR \cong \mC^m \times \mR$ in the same fashion. The m-torus $\mT^m$ is written as
$$\mT^m = \{ \vec{w} = (w_1, \ldots , w_m) \mid w_i \in \mC ~, ~ |w_i| = 1, ~{\rm for ~ all} ~ 1 \leq i \leq m\}$$
Identify $\mT^m$ with a maximal torus in $\SOq$ via its action by coordinate multiplication, for $\vec{z} \in \mC^m$,
\begin{eqnarray}
\vec{w} \cdot \vec{z} & = & (w_1 \cdot z_1, \ldots , w_m \cdot z_m) ~, ~ q = 2m \label{eq-action}\\
\vec{w} \cdot (\vec{z}, x_{2m+1}) & = & (w_1 \cdot z_1, \ldots , w_m \cdot z_m, x_{2m+1}) ~, ~ q = 2m + 1 \nonumber
\end{eqnarray}

Given $\vec{a} = (a_1, \ldots , a_m) \in \mR^m$, set
\begin{equation}\label{eq-eb}
\eb(\vec{a}) = (\exp(2 \pi \sqrt{-1} \cdot a_1), \ldots , \exp( 2 \pi \sqrt{-1} \cdot a_m)) \in \mT^m
\end{equation}

Define a parametrized family of representations, for $0 \leq t \leq 1$,
\begin{equation}
\rho_t^{\alpha} \colon \G \to \mT^m \subset \SOq ~ , ~
\rho_t^{\alpha}(\gamma)(\vec{v}) = \eb(t \cdot \alpha(\gamma)) \cdot \vec{v}  ~ ; \quad \text{set} ~  \rho^{\alpha} = \rho^{\alpha}_1  \label{eq-ortho}
\end{equation}

 \medskip

\begin{prop} Each representation $\alpha \colon \G \to \mR^m$ defines a foliation $\F_{\alpha}$ of $L_0 \times \mS_{\e}^{q-1}$ whose leaves cover $L_0$. Moreover, if the image of $\alpha $ is contained in the rational points $\mQ^m \subset \mR^m$, then all leaves of $\F_{\alpha}$ are compact.
\end{prop}
\sop Given $\alpha \colon \G \to \mR^m$, then the family
$\rho_t^{\alpha}$ is an isotopy from $ \rho^{\alpha}$ to the trivial representation, so $\mS^{q-1}_{\e, \rho} \to L_0$ is bundle-isomorphic to the product bundle $L_0 \times \mS$. If the image of $\alpha$ is contained in $\mQ^m$, then the image of $\rho^{\alpha}$ is a finite subgroup of $\SOq$, hence each leaf of $\F_{\rho}$ is a finite covering of $L_0$ hence is compact. Let $\F_{\alpha}$ denote the foliation which is the image of $\F_{ \rho^{\alpha}}$ under the fiber-preserving diffeomorphism $\mS_{\e, \rho} \cong L_0 \times \mS^{q-1}_{\e}$ induced by the isotopy $\rho_t^{\alpha}$.
\hfill $\eop$

In the case where $q=2$, $\SOq \cong \mS^1$ and $L_0 = \mS^1$, then $\alpha \colon \mZ \to \mR$ is determined by the real number $\alpha = \alpha(1)$, and $\F_{\alpha}$ is the foliation of the 2-torus $\mT^2$ by lines of ``slope'' $\alpha$. Note that for the abstract flat bundle $\mE_{\rho}$ the holonomy is rotation of the fiber $\mS^1$ by the congruence class of $\alpha$ modulo $\mZ$. However, using the isotopy $\rho_t^{\alpha}$ we are able to define the slope of the leaves $\F_{\alpha}$ due to the explicit product structure. A similar phenomenon holds for the general case of $ \rho^{\alpha}$ defined by $\alpha \colon \G \to \mR^m$, although it is not as immediate to imagine the foliation $\F_{\alpha}$ on the product space $L_0 \times \mS^{q-1}_{\e}$. However, there is a standard observation which obviates this problem.

For a vector $\vec{x} \in \mR^k$ let $[\vec{x}] \in \mT^k = \mR^k/ \mZ^k$ denote the coset of $\vec{x} \in \mR^k$. Let $[\vec{0}]$ denote the basepoint defined by the origin in $\mR^k$. The fundamental group $\pi_1(\mT^k, [\vec{0}] )$ acts on the universal cover $\mR^k$ via translations, so is canonically identified with the integer lattice $\mZ^k$.

\begin{prop} \label{prop-realize}
Suppose that $H^1(L_0, \mR) \cong \mR^k$ for $k \geq 1$.
Then there is a surjective map $\pi \colon \G = \pi_1(L_0 , x_0) \to \mZ^k$.
Moreover, there is a smooth map $\tau \colon L_0 \to \mT^k$ with $\tau(x_0) = [\vec{0}] \in \mT^k$, and $\pi = \tau_{*} \colon \pi_1(L_0, x_0) \to \pi_1(\mT^k , [\vec{0}])$.
\end{prop}
\proof
The map $\tau$ is just the period map for $L_0$. More precisely, choose closed $1$-forms $\{\eta_1 , \ldots , \eta_k\}$ on $L_0$ whose cohomology classes are a basis for the image $H^1(L_0, \mZ)  \to H^1(L_0, \mR)$.  Identify the universal cover $\wtL_0$ with the  end-point fixed homotopy classes of paths $\sigma \colon [0,1] \to L_0$   with $\sigma(0) = x_0$.
Each homotopy class of paths admits a smooth representative, so we may assume that  $\sigma$ is smooth. Then define $\wttau \colon \wtL_0 \to \mR^k$ by
$$\wttau(\sigma) ~ = ~  \int_{\sigma} (\eta_1 , \ldots , \eta_k) \, dt ~ \in ~ \mR^k$$
If $\sigma(1) = x_0$, then $\wttau(\sigma) \in \mZ^k$, and so set $\pi([\sigma]) =  \wttau(\sigma)$.
This yields   a well-defined   map $\pi \colon \G \to \mZ^k$, called the period map. Pass  to the quotient  of $\wtL_0$ by the deck translation action of $\G$ to obtain the map  $\tau$. (See \cite{BottTu1982} for details.)
\endproof

\begin{prop} \label{prop-transverse}
Given a foliation $\whF$ of $\mT^k \times \mD^q$ which is transverse to the factor $\mD^q$,
then a smooth map $\tau \colon L_0 \to \mT^k$ induces a foliation $\F$ on $L_0 \times \mD^q$ whose global holonomy map is the composition of $\tau_*$ with the holonomy map $h_{\whF} \colon \mZ^k \to {\rm Diff}(\mD^q)$ of $\whF$.
\end{prop}
\proof
This follows by the standard transversality techniques \cite{CN1985,Milnor2009}.
\endproof

Propositions~\ref{prop-realize} and \ref{prop-transverse} are invoked to obtain a foliation of $M = L_0 \times \mD^q$ given a foliation  of $N_0 = \mT^k \times \mD^q$. Thus, it suffices to construct the examples for this special case $\mT^k \times \mD^q$.
We introduce the basic notations used in the constructions of sequences of such foliations in the following section.

Given $\alpha = (\alpha_1, \ldots, \alpha_m) \colon \mZ^k \to \mR^m$, there is a unique linear extension to $\mR^k$, also denoted by $\alpha \colon \mR^k \to \mR^m$. The extension $\alpha$ is used to define a diagonal action
\begin{eqnarray}\label{eq-affineaction}
T_{\alpha} \colon \mR^k \times (\mT^k \times \mT^m) & \to & (\mT^k \times \mT^m)\\
T_{\alpha}(\vec{\xi})([\vec{x}], [\vec{y}]) & = & ([\vec{x} + \vec{\xi}], [\vec{y} + \alpha(\vec{\xi}) ]) \nonumber
\end{eqnarray}
for $\vec{\xi} \in \mR^k$. The orbits of $T_{\alpha}$ define a foliation $\F_{\alpha}$ on $\mT^k \times \mT^m$ where the leaf through $(x,y) = ([\vec{x}], [\vec{y}]) \in \mT^k \times \mT^m$ is
\begin{equation}
L_{(x,y)} = \{([\vec{x} + \vec{\xi}], [\vec{y} + \alpha(\vec{\xi}) ] )\mid ~ \vec{\xi} \in \mR^k \}
\end{equation}
This is just the standard ``linear foliation'' of $\mT^{k+m}$ by $k$-planes, whose ``slope'' is defined by the map $\alpha$.
If we have a smooth family of representations, $\alpha_t \colon \mZ^k \to \mR^m$, $0 \leq t \leq 1$, then we obtain a smooth family of foliations $\cL_{\alpha_t}$ on $\mT^k \times \mT^m$.

The linear representation $\alpha \colon \mR^k \to \mR^m$ can be composed with the group embedding  of $\mT^m$ into $\SOq$ described above (\ref{eq-action}). We thus  obtain an isometric action $\rho^{\alpha}$ of $\mR^k$ on $\mR^q$, which restricts to actions on each of the subsets
$\mD^q_{\e}$, $\mB^q_{\e}$ and $\mS^{q-1}_{\e}$ of $\mR^q$, and  thus obtain   diagonal actions on the corresponding product spaces. The product actions are also denoted by $T_{\alpha}$. For example, we have
\begin{equation}\label{eq-qaction}
T_{\alpha} \colon \mR^k \times (\mT^k \times \mD^q_{\e}) \to (\mT^k \times \mD^q_{\e})
\end{equation}
The orbits of $T_{\alpha}$ define a foliation of $\mT^k \times \mD^q_{\e}$, again denoted by $\F_{\alpha}$. Note that the action of $T_{\alpha}$ preserves the spherical submanifolds $\mT^k \times \mS^{q-1}_{\e'}$ for $0 \leq \e' \leq \e$.

\section{Constructing the plug} \label{sec-plug}

In this section, we inductively define a sequence of foliations $\whF_{\ell}$ of $\mT^k \times \mD^{q}$, for $\ell \geq 0$. In Section~\ref{sec-estimates}, we give criteria for this sequence to converge to a $C^r$ foliation $\whF$.

We first  establish some   notations.
Let $\vec{e}_{j} = (0, \ldots, 1, \ldots, 0) \in \mZ^k$ be the standard basis vector, where the sole non-zero entry is in the $j^{th}$-coordinate. Then   $\{\vec{e}_1, \ldots , \vec{e}_k\}$  is called the \emph{standard basis}.

Let $\alpha_0$ denote the trivial representation, with $\alpha_0(\vec{e}_{j}) = \vec{0}$ for all $1 \leq j \leq k$. Then
for each $\ell \geq 1$, choose a representation $\alpha_{\ell} \colon \mZ^k \to \mQ^m$. Let $\whalpha = \{\alpha_{\ell} \mid \ell \geq 0\}$ denote the sequence of representations.

 Set $\alpha_{\ell}(\vec{e}_{j}) =\vec{a}_{\ell, j} = [a_{\ell, j,1} , \ldots , a_{\ell, j,m}]^t \in \mQ^m$.

 We let $\tn \cdot \tn$ denote the sup-norm of any collection of real numbers. For example, we have
\begin{equation}\label{eq-tn}
\tn \alpha_{\ell} \tn ~ = ~ \max \{ |a_{\ell, j,m}| \, \mid 1 \leq j \leq k , 1 \leq k \leq m\}
\end{equation}

In this case,  the quantity $\tn \alpha_{\ell} \tn$ is a measure of the ``maximum slope'' of the flat bundle foliation obtained from the representation $\rho_{\ell} \equiv  \rho^{\alpha_{\ell}} = \eb(\alpha_{\ell})$.

Define representations $\rho_{\ell} = \rho^{\alpha_{\ell}} \colon \mZ^k \to \mT^m \subset \SOq$   as in (\ref{eq-ortho}),  setting
\begin{equation}\label{eq-defrho}
\rho_{\ell}(\gamma)(\vec{v}) = \eb(\alpha_{\ell}(\gamma)) \cdot \vec{v} \quad , ~ \text{for} ~ \gamma \in \mZ^k ~, ~ \vec{v} \in \mR^q
\end{equation}
Let $\Delta_{\ell} \subset \mQ^m$ denote the image of $\alpha_{\ell}$.
Then the  image of $\rho_{\ell}$ is isomorphic to the finite subgroup, $\fG_{\ell} \equiv  \Delta_{\ell}/\mZ^m \subset \mQ^m/\mZ^m$. Let $d_{\ell} = \#\fG_{\ell}$ denote the order of $\fG_{\ell}$.
Note that   $\fG_0$ is the trivial group, as the image of $\alpha_0$ is all of $\mZ^k$.  We assume that all other subgroups $\fG_{\ell}$ for $\ell \geq 1$ are non-trivial.

Let $\Lambda_{\ell} = \ker \{\rho_{\ell} \colon \mZ^k \to \SOq\} \subset \mZ^k$. Then $\fG_{\ell} \cong \mZ^k/\Lambda_{\ell}$, where   $\Lambda_{\ell}$ is free abelian with rank $k$.

For each $\ell \geq 1$, choose a set of generators $\{e^{\ell}_1, \ldots , e^{\ell}_k \} \subset \Lambda_{\ell}$ with positive orientation for
$ \Lambda_{\ell}$, which determine an isomorphism $\phi_{\ell} \colon \mZ^k \to \Lambda_{\ell}$ given by $\phi_{\ell}(c_1, \ldots , c_k) = \sum_{j=1}^{k} c_{j} e^{\ell}_{j}$.
Let $\phi_{\ell} \colon \mR^k \to \mR^k$ also denote the extension of $\phi_{\ell}$ to a linear isomorphism of $\mR^k$.  The
 induced isomorphism on quotients is denoted by $\pb_{\ell} \colon \mT^k = \mR^k/\mZ^k \to \mR^k/\Lambda_{\ell}$.

The collection of subgroups $\Lambda_{\ell}  \subset \mZ^k$ gives rise to a descending chain as follows.
For $\ell \geq 1$, define
\begin{equation}\label{eq-Phi}
 \Phi_{\ell} \equiv \phi_1 \circ \cdots \circ \phi_{\ell} \colon \mZ^k \to \mZ^k \quad ; ~ \G_{\ell} = {\rm Image}\{\Phi_{\ell}\} \subset \mZ^k
\end{equation}
Set $\G_0 = \mZ^k$. Then for  each $\ell \geq 0$,  $\G_{\ell +1} \subset \G_{\ell}$.
Note that  $\G_1 = \Lambda_1$. For $\ell \geq 1$, we have
the quotient $\G_{\ell -1}/\G_{\ell} \cong \Lambda_{\ell }$.
 It follows that     $\G_{\ell}$ is  a subgroup of finite index,   $q_{\ell} \equiv [\G_{\ell} : \mZ^k] = d_1 \cdots d_{\ell}$.

 Set $M_{\ell} = \mR^k/\G_{\ell}$, and let $p_{\ell} \colon M_{\ell} \to M_{\ell -1}$ be the quotient map, which is a locally isometric  covering map of degree $d_{\ell}$ with Galois group $\Lambda_{\ell}$.
The composition of covering maps,
 $\pi_{\ell} \colon \mR^k/\G_{\ell} \to \mR^k/\G_0$, has Galois group $\mathfrak{k}_{\ell} = \G_0 /\G_{\ell}$ which has index $q_{\ell}$.

 \medskip

The initial step in our construction is to     realize a  given representation  $\rho_{\ell}$  as the holonomy restricted to an invariant disk bundle, for  a foliation  on $\mT^k \times \mD^q$ which is a product near the boundary. For this, we introduce the non-linear aspect of the construction.
Choose a non-increasing smooth function $\mu \colon [0,1] \to [0,1]$ such that
\begin{equation}\label{eq-mudef}
\mu (s) = \begin{cases}
1 & \text{if } ~ 0 \leq s \leq 2/3 \\
0 & \text{if } ~ 3/4 \leq s \leq 1
\end{cases}
\end{equation}
For each $\ell \geq 1$, we extend the representation $\alpha_{\ell} \colon \mZ^k \to \mQ^m$ to a continuous family
$\alpha_{\ell,t} \colon \mZ^k \to \mR^m$, $0 \leq t \leq 1$, by setting $\alpha_{\ell,t}(\vec{e}_{j}) = \mu(t) \cdot \vec{a}_{\ell,j}$ for $1 \leq j \leq k$.
Define a   smooth family of representations,
\begin{equation}\label{eq-famrep}
\rho_{\ell,t} \colon \mZ^k \to \mT^m ~, \quad \rho_{\ell,t,j}(\vec{v}) = \rho_{\ell,t}(\vec{e}_j)(\vec{v}) = \eb(\mu(t) \cdot \vec{a}_{\ell,j}) \cdot \vec{v}
\end{equation}
Note that for $3/4 \leq t \leq 1$, $\rho_{\ell,t}$ is the identity $I \in \SOq$, and for $0 \leq t \leq 2/3$ the action of $\rho_{\ell,t,j}$ is multiplication by $\eb(\vec{a}_{\ell,j})$, hence $\rho_{\ell,t}$ has image isomorphic to $\fG_{\ell}$.

A point $\vec{v} \in \mR^q$ is \emph{generic} for $\rho_{\ell}$ if $\rho_{\ell}(\gamma) \cdot \vec{v} \ne \vec{v}$ for all $\gamma \not\in \Lambda_{\ell}$.
Let $\cO_{\ell} \subset \mR^q$ denote the generic points for $\rho_{\ell}$.
The fix-point set for an isometry $\rho_{\ell}(\gamma)$ of $\mR^q$ is a proper subspace of $\mR^q$ if $\gamma \not\in \Lambda_{\ell}$, and as $\fG_{\ell}$ is a finite group, the set of non-generic points in $\mR^q$ is a finite union of proper subspaces. Thus, $\cO_{\ell}$ is an open and dense subset.

We now begin the inductive construction of foliations $\{\whF_{\ell} \mid \ell \geq 0\}$ on $\mT^k \times \mD^q$. Define $\whF_0$ to be  the product foliation   on $\mT^k \times \mD^{q}$.
We   present the  construction of $\whF_1$ from $\whF_0$ before giving the general inductive step, as this case is technically simpler than the general case.

It is helpful to keep in mind that the following   constructions of maps generalizes the standard construction for dyadic solenoids, and  the algebraic steps we describe all have geometric interpretations, albeit in higher dimensions. The figure below is the standard
illustration of the imbedding of the solid $2$-torus $\mT^2$ into itself; or in the language to follow, this represents the map $\Psi_1 \colon N_1 \to N_0$.

\begin{figure}[H]
\includegraphics[width=.4\textwidth]{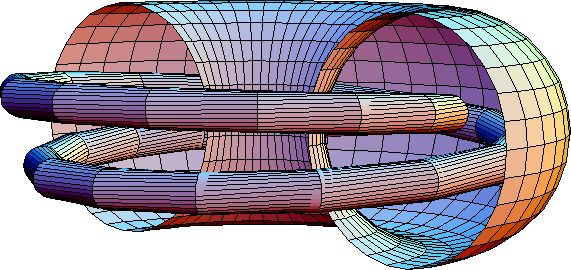}
\caption{Embedded solid torus}
\label{figure-originalhirsch}
\vspace{-20pt}
\end{figure}

For $\ell = 0$, set $\e_0 = 1$,   $\e_0' = 2/3 \cdot \e_0$ and $\e_0'' = 3/4 \cdot \e_0$. Then set
$$  \rho_0 \equiv Id ~ , ~ \Phi_0 \equiv Id ~, ~ \G_0 \equiv \Lambda_0 = \mZ^k ~ , ~ N_0 \equiv \mT^k \times \mD^{q}_{\e_0} = \mT^k \times \mD^{q}_{1} ~ , ~ \mK_0 \equiv \mD^{q}_{\e_0} = \mD^{q}_{1}$$
For step $\ell =0$, we   make the following  formal definitions. Set $\whN_0 = N_0$ and let $\Psi_0 \colon N_0 \to \whN_0$ be the identity map. The foliation $\F_0$ of $N_0$ is the product foliation, and   define $\whF_0$ on $\whN_0$ to be  the image of $\F_0$ under $\Psi_0$.
Then  $ \whL_0 = \Psi_0(\mT^k \times \vec{0}) \subset \whN_0$ is the leaf of $\whF_0$ through $\vec{z}_0 = \vec{0}$.
Set $\whmK_0 = \Psi_0([\vec{0}] \times K_0)$.

For $\ell =1$, we next construct the foliation $\whF_1$ of $\whN_0$.
To begin, there is a continuous decomposition
\begin{equation}\label{eq-sphere0}
N_0 = \bigcup_{0 \leq r \leq \e_0} ~ \mT^k \times \mS^{q-1}_{r}
\end{equation}
For each $0 \leq r \leq \e_0$, let $\F_1$ restricted to $\mT^k \times \mS^{q-1}_{r}$ be the foliation defined by the representation $\rho_{1,t}$ where $t = r/\e_0$.
The foliation $\F_1$ of $N_0$ is smooth, as $\rho_{1,t}$ depends smoothly on $t$.

The family $\{\rho_{1,t} \mid 0 \leq t \leq 1\}$ is an isotopy between $\rho_1$
and the trivial representation $\rho_0$.
The foliation $\F_1$ restricted to $\mT^k \times (\mD^{q}_{\e_0} - \mB^{q}_{\e_0''})$ is the product foliation; $\F_1$ restricted to $\mT^k \times \mD^{q}_{\e_0'}$ equals $\F_{\rho_1}$ whose holonomy is given by multiplication by the complex vectors $\rho_1(\gamma) \in \mT^m$, for $\gamma \in \mZ^k$. For $\e_0' < r < \e_0''$, the foliation restricted to $ \mT^k \times \mS^{q-1}_{r}$ is the suspension of an isometric action.

Note that $\F_1$ is a distal foliation. Given any two points $z \ne z' \in \mS^{q-1}_{r}$ the holonomy action of $\F_1$ on these points is isometric, hence stay  a bounded distance apart. On the other hand, if $z \in \mS^{q-1}_{r}$ and $z' \in \mS^{q-1}_{r'}$ for $0 \leq r < r' \leq \e_0$ then their orbits remain on distinct spherical shells $\mS^{q-1}_{r}$ and $\mS^{q-1}_{r'}$, hence remain a bounded distance apart.

We have now constructed the foliation $\F_1$ on $N_0$, which is the ``standard solid
torus''.   Let $\whF_1$ be the foliation of $\whN_0$ which is
the image of $\F_1$ under $\Psi_0 \colon N_0 \to \whN_0$. As
$\Psi_0$ is the identity, this step is again purely formal.

The inductive step for $\ell=1$ includes several further choices of data.

Let $\vec{z}_1 \in \mS^{q-1}_{\e_0/2} \cap \cO_1 \subset \mB^{q}_{\e_0'} \cap \cO_1$ be a generic point for $\rho_1$.
For $\gamma \in \mZ^k$ set $\vec{z}_{1,\gamma} = \rho_1(\gamma) (\vec{z}_1)$.
The $\rho_1$-orbit $\{ \vec{z}_{1,\gamma} \mid \gamma \in \mZ^k \}$ of $\vec{z}_1$ is finite with order $d_1$, so there exists $\e_1 > 0$ such that the closed disk centered at $\vec{z}_1$ satisfies
\begin{equation}\label{eq-disk2}
\mD^q_{\e_1}(\vec{z}_1) \equiv \{ \vec{z} \in \mR^q \mid \|\vec{z} - \vec{z}_1 \| \leq \e_1 \} \subset \mB^{q}_{\e_0'} \cap \cO_1
\end{equation}
and the translates under the action of $\rho_1$ are disjoint. That is, if $\gamma \in \mZ^k$ satisfies $\gamma \not\in \Lambda_1$ then
$$\mD^q_{\e_1}(\vec{z}_1) \cap \rho_1(\gamma) \cdot \mD^q_{\e_1}(\vec{z}_1) ~ = ~ \mD^q_{\e_1}(\vec{z}_1) \cap \mD^q_{\e_1}(\vec{z}_{1, \gamma}) ~ = ~ \emptyset$$

Note that $\|\vec{z}_1\| = \e_0/2$ and $\ds \mD^q_{\e_1}(\vec{z}_1) \subset \mB^{q}_{\e_0'} \cap \cO_1$ implies that $\e_1 < \e_0/6$ and hence $\mD^q_{\e_0/3} \cap \mD^q_{\e_1}(\vec{z}_1) = \emptyset$.
The finite union of the translates of $\mD^q_{\e_1}(\vec{z}_1)$ under the action $\rho_1$ is denoted by
\begin{equation}\label{eq-K1}
K_1 = \bigcup_{\gamma \in \mZ^k} ~ \mD^q_{\e_1}(\vec{z}_{1,\gamma}) ~ \subset ~ K_0 = \mD^{q}_{1}
\end{equation}
Here is an illustration of the set $K_1$ in the case  $q=2$ and $\fG_1 = \mZ/6\mZ$ the  cyclic group  of order $6$.
The points $z_{1,i}$ label the orbit of $\vec{z}_1$ under $\rho_1$.

\begin{figure}[H]
\includegraphics[width=0.28\textwidth]{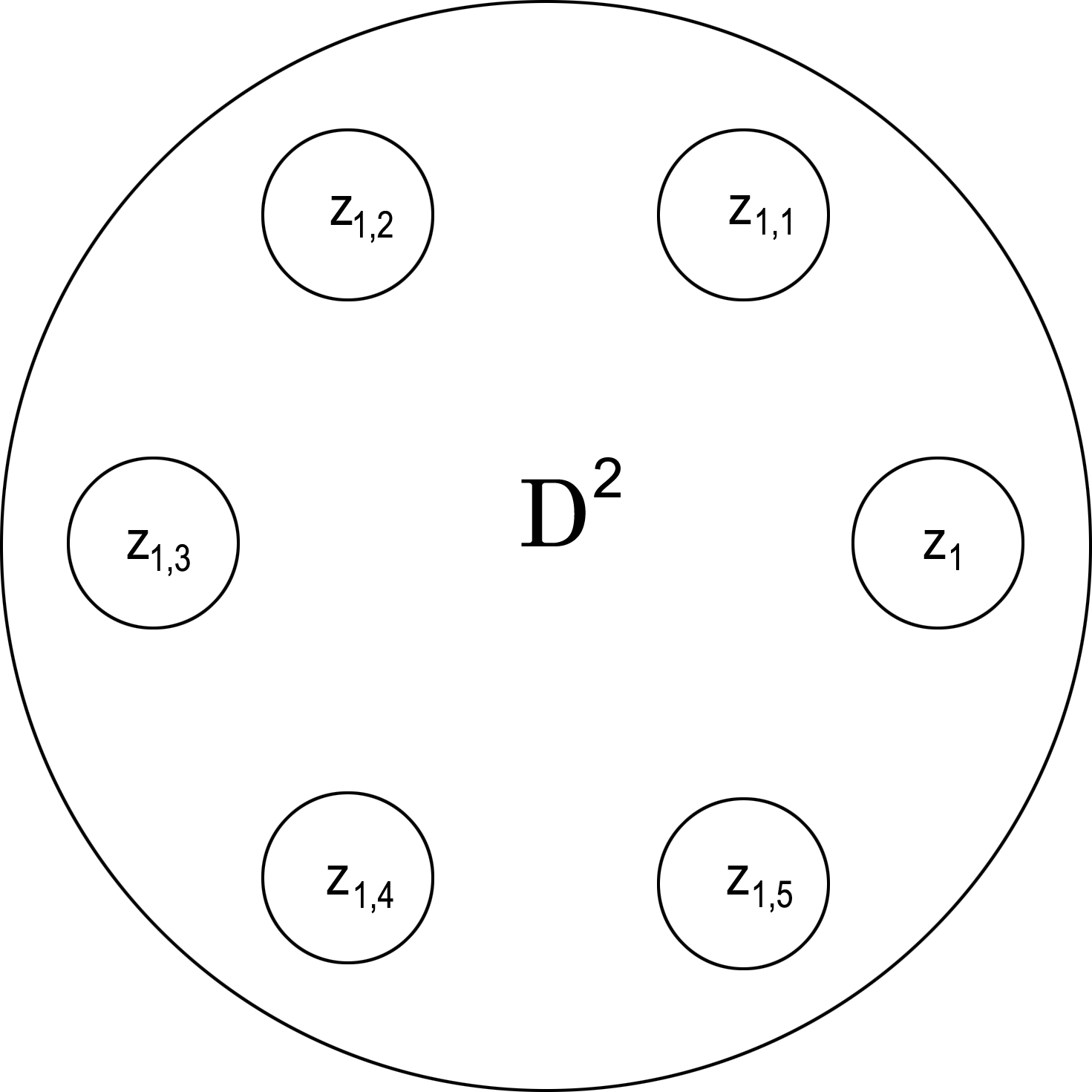}
\caption{Invariant family of disks}
\label{figure-disks}
\vspace{-20pt}
\end{figure}

Let $N_1'$ denote the $T_{\alpha_1}$-saturation of the set $[\vec{0}] \times \mD^q_{\e_1}(\vec{z}_1) \subset N_0$. That is,
\begin{equation}
N_1' = \bigcup_{\vec{\xi} \in \mR^k} ~   [\vec{\xi}] \times \rho_1(\vec{\xi} ) \cdot \mD^q_{\e_1}(\vec{z}_1)  ~ \subset ~ N_0 = \mT^k \times \mD^{q}_{1}
\end{equation}
Finally, we put the foliation $\F_1$ of $N_1'$ into ``standard form'' in preparation for the next stage of the induction.
Define a $T_{\alpha_1}$-equivariant map, for $\xi \in \mR^k$
\begin{eqnarray}
\varphi_1 \colon \mR^k/\Lambda_1 \times \mD^q_{\e_1} & \longrightarrow & N_1'   \nonumber \\
\varphi_1([\vec{\xi}], \vec{z}) & = &   [\vec{\xi}] \times \rho_1(\vec{\xi})( \vec{z} + \vec{z}_1)  \label{eq-holo1}
\end{eqnarray}
This is well-defined since  $\Lambda_1 = \ker \rho_1$.
Note that the product foliation on $\mR^k/\Lambda_1 \times \mD^q_{\e_1} $ is mapped by $\varphi_1$ to the restriction of $\F_1$ to $N_1'$.

Recall that  $\phi_{\ell} \colon \mR^k \to \mR^k$  is the extension of the homomorphism $\phi_{\ell} \colon \mZ^k \to \Lambda_{\ell}$ defined by the choice of basis for $\Lambda_{\ell}$, and  $\pb_{\ell} \colon \mT^k = \mR^k/\mZ^k \to \mR^k/\Lambda_{\ell}$ is the induced map on quotients.
For $\ell =1$, extend  this to a diffeomorphism
\begin{equation}\label{eq-normalize}
\phi_1 \colon N_1 \equiv \mT^k \times \mD^q_{\e_1}   \to   \mR^k/\Lambda_1 \times \mD^q_{\e_1} ~ , \quad ([\vec{x}],  \vec{y} )   \mapsto   (\phi_1[\vec{x}] ,  \vec{y} )
\end{equation}
Let $\psi_1 = \varphi_1 \circ \phi_1 \colon N_1 \to N_1'$, and set  $\Psi_1 = \Psi_0 \circ \psi_1 \colon N_1 \to \whN_0$ which maps the product foliation on $N_1$  to the restriction of  $\whF_1$ to  $\whN_1$. We have the diagram:
\begin{equation}\label{eq-Psi1}
\Psi_1 \colon N_1 \equiv \mT^k \times \mD^q_{\e_1}    \stackrel{\phi_1}{\longrightarrow}   \mR^k/\Lambda_1 \times \mD^q_{\e_1}  \stackrel{\varphi_1}{\longrightarrow}   N_1'  \subset   N_0 = \mT^k \times \mD^{q}_{\e_0} \stackrel{\Psi_0}{\longrightarrow}  \whN_0
\end{equation}
Let $\whN_1$ denote the image of $\Psi_1$ and set
 $\whmK_1 = \whN_1 \cap \whmK_0$, which is the  union of $d_1$ closed disks.

Let $L_1 \subset N_1'$ be the leaf which is the $T_{\alpha_1}$-orbit of $[\vec{0}] \times \vec{z}_1$, then $\whL_1 = \Psi_0(L_1) = \Psi_1(\mT^k \times \vec{0})$  is a leaf of $\whF_1$ and $\whN_1$ is the closed $\e_1$-disk bundle about $\whL_1$.
This concludes   step $\ell =1$.

\medskip

Now assume that the foliation $\whF_{\ell}$ of $\whN_0$ has been constructed,  and we construct the foliation $\whF_{\ell +1}$.

Assume that points
$\{\vec{z}_1, \ldots , \vec{z}_{\ell}\} \subset \mD^q_{\e_0}$ have been chosen, as well as the sequence
$0 < \e_{\ell} < \cdots < \e_1 < 1$ where $\e_{k+1} < \e_k/6$ for $1 \leq k < \ell$.
Let $\e_k' = 2/3 \cdot \e_k$ and $\e_k'' = 3/4 \cdot \e_k$ for $1 \leq k \leq \ell$.

By the inductive hypotheses, we are also given $\Psi_{\ell} = \Psi_{\ell -1} \circ \psi_{\ell} $ such that
\begin{equation}\label{eq-Psiell}
\Psi_{\ell}  \colon N_{\ell} \equiv \mT^k \times \mD^q_{\e_{\ell}}    \stackrel{\phi_{\ell}}{\longrightarrow}   \mR^k/\Lambda_{\ell} \times \mD^q_{\e_{\ell}}  \stackrel{\varphi_{\ell}}{\longrightarrow}   N_{\ell}'  \subset   N_{\ell -1} = \mT^k \times \mD^{q}_{\e_{\ell -1}} \stackrel{\Psi_{\ell -1}}{\longrightarrow}  \whN_0
\end{equation}
where  $\whN_{\ell}$ denotes the image of $\Psi_{\ell}$ and
 $\whmK_{\ell} = \whN_{\ell} \cap \whmK_0$, which is the  union of $q_{\ell}$ closed disks.

The foliation $\whF_{\ell}$ of $\whN_0$ restricted to $\whN_{\ell}$ is the image under $\Psi_{\ell}$ of the product foliation on $N_{\ell}$.

The composition (\ref{eq-Psiell}) realizes an important technical aspect of our construction. The image of the core leaf $\mT^k \times \{0\}$ of the product foliation
on $N_{\ell}$ is mapped by  $\Psi_{\ell}$ to a leaf of $\whF_{\ell}$. Moreover,  $\Psi_{\ell}$ also maps the trivial framing of the disk bundle to a twisted framing (or synchronous) framing of a tubular neighborhood of the image leaf. This allows us to make the ``foliated surgery'' described below on a product bundle, which reduces to a problem exactly analogous to the case $\ell =1$.

Topologically, the manifold with boundary $\whN_{\ell} \subset \whN_0$ is the $\ell^{th}$-iteration of an embedding of a solid torus in the initial solid torus $\whN_0$. For $k=1$, this is visualized as the iteration of the embedding in   Figure 1.
For $k > 1$, the embedding depends upon the sequence of vectors
 $\alpha_{i} \in \mathrm{Hom}(\mZ^k ,  \mR^m)$ for $1 \leq i \leq \ell$, and is essentially impossible to visualize.

The first step in the construction of $\whF_{\ell + 1}$ is to construct the foliation $\F_{\ell + 1}$ on $N_{\ell}$ with the property that its restriction to  an open neighborhood of $\partial N_{\ell}$ is the product foliation.
We proceed as for the case $\ell =1$. There is a continuous decomposition
\begin{equation}
N_{\ell} = \bigcup_{0 \leq r \leq \e_{\ell}} ~ \mT^k \times \mS^{q-1}_{r}
\end{equation}
For each $0 \leq r \leq \e_{\ell}$, let $\F_{\ell+1}$ restricted to $\mT^k \times \mS^{q-1}_{r}$ be the foliation defined by the representation $\rho_{\ell +1,t}$ where $t= r/\e_{\ell}$.
This defines the smooth  foliation $\F_{\ell +1}$ of $N_{\ell}$ and we check its properties.

The family $\{\rho_{\ell+1 , t} \mid 0 \leq t \leq 1\}$ is an isotopy between $\rho_{\ell + 1}$
and the trivial representation $\rho_0$.
Thus, the foliation $\F_{\ell + 1}$ restricted to $\mT^k \times (\mD^{q}_{\e_{\ell}} - \mB^{q}_{\e_{\ell}''})$ is the product foliation.

The restriction of  $\F_{\ell + 1}$   to $\mT^k \times \mD^{q}_{\e_{\ell}'}$ is the foliation $\F_{\rho_{\ell + 1}}$,
whose holonomy action on  $\ds [\vec{0}] \times \mD^{q}_{\e_{\ell}'}$
is given by multiplication by the complex vectors $\rho_{\ell + 1}(\gamma) \in \mT^m$, for $\gamma \in \mZ^k$.

For $\e_{\ell}' < r < \e_{\ell}''$, the foliation restricted to $ \mT^k \times \mS^{q-1}_{r}$ is the suspension of an isometric action.
  Thus, as before, $\F_{\ell + 1}$ is a distal foliation of $N_{\ell}$.

Define the foliation $\whF_{\ell +1}$ on $\whN_0$   to be $\whF_{\ell}$ on $\whN_0 - \whN_{\ell}$, and  $\whF_{\ell +1}$ on   $\whN_{\ell}$   is  the image of $\F_{\ell + 1}$.

By the inductive hypothesis, the image under $\Psi_{\ell}$ of the product foliation on $N_{\ell}$ equals the restriction of $\whF_{\ell}$ on $\whN_{\ell}$, hence the image under $\Psi_{\ell}$ of $\F_{\ell +1 }$ on $\whN_{\ell}$ agrees with $\whF_{\ell}$ on an open neighborhood of $\partial \whN_{\ell}$, and so $\whF_{\ell +1}$  is a smooth   foliation of $\whN_0$.

 It remains to set up the remaining data for the induction.

Let $\vec{z}_{\ell + 1} \in \mS^{q}_{\e_{\ell}/2} \subset \mB^{q}_{\e_{\ell}'} \cap \cO_{\ell + 1}$ be a generic point for $\rho_{\ell +1}$.
Then set
\begin{equation}\label{eq-defzell}
\vec{z}_{\ell +1,\gamma} = \rho_{\ell+1}(\gamma) (\vec{z}_{\ell + 1}) ~ , \quad \gamma \in \mZ^k
\end{equation}
The $\rho_{\ell + 1}$-orbit $\{ \vec{z}_{\ell + 1,\gamma} \mid \gamma \in \mZ^k \}$ of $\vec{z}_{\ell + 1}$ is finite, so there exists $\e_{\ell + 1} > 0$ such that the closed disk centered at $\vec{z}_{\ell + 1}$ satisfies
\begin{equation}\label{eq-diskell}
\mD^q_{\e_{\ell + 1}}(\vec{z}_{\ell + 1}) \subset \mB^{q}_{\e_{\ell}'} \cap \cO_{\ell + 1}
\end{equation}
and the translates under the action of $\rho_{\ell + 1}$ are disjoint.
Note that $\ds \mD^q_{\e_{\ell + 1}}(\vec{z}_{\ell + 1}) \subset \mB^{q}_{\e_{\ell}'} \cap \cO_{\ell + 1}$ implies that $\e_{\ell + 1} < \e_{\ell}/6$ and hence $\mD^q_{\e_{\ell}/3} \cap \mD^q_{\e_{\ell + 1}}(\vec{z}_{\ell + 1}) = \emptyset$.
The finite union of the translates of $\mD^q_{\e_{\ell + 1}}(\vec{z}_{\ell + 1})$ under the action $\rho_{\ell + 1}$ is denoted by
\begin{equation}\label{eq-Kell+1}
K_{\ell + 1} = \bigcup_{\gamma \in \mZ^k} ~ \mD^q_{\e_{\ell + 1}}(\vec{z}_{\ell + 1, \gamma})
\end{equation}
 Then $K_{\ell + 1}$ is the disjoint union of $d_{\ell +1}$ closed disks, each of radius $\e_{\ell +1}$, so that   $K_{\ell+1} \cap \mD^q_{\e_{\ell}/3} = \emptyset$.

Let $N_{\ell + 1}'$ denote the $T_{\alpha_{\ell + 1}}$-saturation of the set $[\vec{0}] \times \mD^q_{\e_{\ell + 1}}(\vec{z}_{\ell + 1}) \subset N_{\ell}$. That is,
\begin{equation}
N_{\ell + 1}' = \bigcup_{\vec{\xi} \in \mR^k} ~  [\vec{\xi}] \times  \rho_{\ell + 1}(\vec{\xi} ) \cdot \mD^q_{\e_{\ell + 1}}(\vec{z}_{\ell + 1})  ~ \subset ~ N_{\ell}
\end{equation}
Define a $T_{\alpha_{\ell + 1}}$-equivariant map, for $\xi \in \mR^k$
\begin{eqnarray}
\varphi_{\ell + 1} \colon \mR^k/\Lambda_{\ell + 1} \times \mD^q_{\e_{\ell + 1}} & \longrightarrow & N_{\ell + 1}' \label{eq-Fell+1}\\
\varphi_{\ell + 1}([\vec{\xi}], \vec{z}) & = & ( [\vec{\xi}], \rho_{\ell + 1}(\vec{\xi})(\vec{z}_{\ell + 1} + \vec{z}) ) \nonumber
\end{eqnarray}
This is well-defined precisely because $\Lambda_{\ell + 1} = \ker \rho_{\ell + 1}$.
The product foliation on $\mR^k/\Lambda_{\ell + 1} \times \mD^q_{\e_{\ell + 1}} $ is mapped by $\varphi_{\ell + 1}$ to the restriction of $\F_{\ell + 1}$ to $N_{\ell + 1}'$.

Extend the map $\pb_{\ell} \colon \mT^k = \mR^k/\mZ^k \to \mR^k/\Lambda_{\ell}$ above
to a diffeomorphism
\begin{equation} \label{eq-normalizeell}
\phi_{\ell + 1} \colon  N_{\ell + 1} \equiv  \mT^k \times \mD^q_{\e_{\ell + 1}}  \to  \mR^k/\Lambda_{\ell + 1} \times \mD^q_{\e_{\ell + 1}}   \quad , ~  ([\vec{x}], [\vec{y}])   \mapsto   (\phi_{\ell + 1}[\vec{x}] , [\vec{y}])
\end{equation}
Let $\psi_{\ell + 1} = \varphi_{\ell + 1} \circ \phi_{\ell + 1} \colon N_{\ell + 1} \to N_{\ell + 1}' \subset N_{\ell}$.
Set $\Psi_{\ell + 1} = \Psi_{\ell} \circ \psi_{\ell + 1} \colon N_{\ell + 1} \to \whN_0$. The image of $\Psi_{\ell + 1}$ is denoted by $\whN_{\ell + 1}$ which is a closed subset of $\whN_0$.
Then $\Psi_{\ell + 1}$ maps to product foliation on $N_{\ell + 1}$ to the restriction of $\whF_{\ell + 1}$ to $\whN_{\ell + 1}$.
Set $\whmK_{\ell + 1} = \whN_{\ell + 1} \cap \whmK_{\ell} = \whN_{\ell + 1} \cap \whmK_{0}$ which is a union of $q_{\ell +1} = d_1 \cdots d_{\ell +1}$ closed disks, each with   radius $\e_{\ell +1}$.

Let $L_{\ell + 1} \subset N_{\ell + 1}'$ be the leaf of $\F_{\ell + 1}$ given by the
$T_{\alpha_{\ell + 1}}$-orbit of $[\vec{0}] \times \vec{z}_{\ell + 1}$,
and set $\whL_{\ell + 1} = \Psi_{\ell}(L_{\ell + 1}) = \Psi_{\ell + 1}(\mT^k \times \vec{0})$.
Then $\whL_{\ell + 1}$ is a leaf of $\whF_{\ell + 1}$
and $\whN_{\ell + 1}$ is a closed $\e_{\ell + 1}$-disk bundle about $\whL_{\ell + 1}$.

This completes the induction. Note that we obtain as a result:
\begin{enumerate}
\item a sequence of nested compact $(k+q)$-dimension submanifolds with boundary,
$$\mT^k \times \mD^q_{\e_0} \equiv \whN_0 \supset \whN_1 \supset \cdots \supset \whN_{\ell} \supset \cdots $$
\item a sequence of nested compact $q$-dimension submanifolds with boundary, $\whmK_{\ell} = \whN_{\ell} \cap \whmK_0$,
$$ \whmK_0 \supset \whmK_1 \supset \cdots \supset \whmK_{\ell} \supset \cdots $$
Moreover, $\whmK_{\ell}$ is a union of $q_{\ell}$ closed  disks, each with radius $\e_{\ell} \leq  6^{-\ell}$. \medskip
\item a sequence of smooth foliations $\whF_{\ell}$ of $\whN_0$ such that $\whF_{\ell'} = \whF_{\ell} \mid (\whN_0 - \whN_{\ell})$ for all $\ell' > \ell$.
\end{enumerate}

The intersection $ \cS = \bigcap_{\ell \geq 0} ~ \whN_{\ell} $ is homeomorphic to a   solenoid as defined in Definition~\ref{def-solenoid}, and
the intersection
$ \mK_* = \bigcap_{\ell \geq 0} ~ \whmK_{\ell} = \cS \cap \whmK_0$
is a Cantor set.

Note that at each stage of the induction, the restriction of $\whF_{\ell}$ to the embedded torus
$\ds \Psi_{\ell} \left( \mT^k \times \mD^q_{\e_{\ell}/3}  \right)$
is a foliation with all leaves compact.

\section{$C^r$-norms}\label{sec-norms}

For a vector-valued  function $f$ defined on an open subset $U \subset \mR^q$, let $\nabla_j f = \nabla_{\vec{e}_j} f$ denote the partial derivative in the direction of the basis vector $\vec{e}_j$ for $1 \leq j \leq q$. Given a string
$J = (j_1, j_2, \ldots , j_{r})$ with values $1 \leq j_k \leq q$, denote the corresponding partial derivative of order $|J| \equiv r$ by
$$\nabla_J f = \nabla_{j_{r}} \circ \cdots \circ \nabla_{j_{1}} f$$
If $r = 0$,  so $J$ is the empty string, then $\nabla_J$ is just the identity map.

Introduce  the uniform $C^r$-semi-norms $\| f \|_r'$ and norm $\| f \|_r$ on $U$, defined by
\begin{equation}
\| f \|_r' = \sup_{\vec{x} \in U} ~ \left\{ \max_{|J| = r} ~  \| \nabla_J f (\vec{x}) \| \right\} ~, \quad
\| f \|_r = \sup_{\vec{x} \in U} ~ \left\{ \max_{|J| \leq r} ~  \| \nabla_J f (\vec{x}) \| \right\}
\end{equation}

\medskip

The   construction in section~\ref{sec-plug}    are based on a ``standard deformation'':
\begin{equation}\label{eq-hr}
\text{for}~ \vec{a} \in \mR^m ~, \quad g_{\vec{a}} (\vec{z}) = \begin{cases}
\eb(\vec{a}) \cdot \vec{z} & \text{if } ~ 0 \leq t \leq 2/3\\
\eb(\mu(t) \cdot \vec{a}) \cdot \vec{z} & \text{if } ~ 2/3 < t < 3/4\\
\vec{z} & \text{if } ~ 3/4 \leq t \leq 1
\end{cases} ~, ~ t = \|\vec{z}\|
\end{equation}
Note that $g_{\vec{a}}$ is the identity outside of $\mD^q_{3/4}$, and is the constant ``rotation''  by $\eb(\vec{a})$ on $\mD^q_{2/3}$.
The function $\mu \colon [0,1] \to [0,1]$ is assumed to be smooth, and vanishes for $t \leq 2/3$, so the composition
$\whmu(\vec{z})  \equiv \mu(\|\vec{z}\|)$ is smooth on the compact disk $\mD^q$. Introduce constants
$B_k ~ = ~ \| \whmu \, \|_k'$ for all $k \geq 0$.

\medskip

The following is the key technical estimate, which along with appropriate rescaling arguments, is used to estimate   the $C^r$-norms of the  sequence of foliations $\{\whF_{\ell} \mid \ell \geq 0\}$.

\begin{lemma}\label{lem-muest}
For all integers $r \geq 0$,
there exists $C_{r} > 0$ such that for   $\tn  \vec{a} \tn \leq 1$,
\begin{equation}\label{eq-normest}
 \| g_{\vec{a}} - Id \|_{r}' \leq     C_{r} \cdot  \tn \vec{a} \tn
\end{equation}
For  the special cases  $C_0 = 2\pi$  and $C_1 = 2\pi (1 + B_1)$, there is no restriction on  $\tn  \vec{a} \tn$. In general set,
\begin{equation}\label{eq-normmax}
\whC_{\ell} \equiv \max \{C_0, C_1, \ldots, C_{\ell }\}
\end{equation}
\end{lemma}
\proof
Note that for $\|\vec{z}\| > 3/4$ we have $ g_{\vec{a}} - Id = 0$, so the result is trivial in this case.

For $\|\vec{z}\| = t \leq  3/4$, using  the observation  that  for all $s \in \mR$, $\|\exp (2\pi \, s)  - \exp (0) \| \leq 2\pi \, |s|$, then
$$
\|g_{\vec{a}}(\vec{z}) - \vec{z} \|   =    \| \eb(\mu(t) \cdot \vec{a}) \cdot \vec{z} -\vec{z} \, \|
  \leq   2\pi \tn \vec{a}\tn \cdot   \| \vec{z} \|  \leq   2 \pi \tn \vec{a} \tn
$$

Thus, $C_0 = 2\pi$ satisfies (\ref{eq-normest}) for $r =0$.
 Next, observe that    for $1 \leq j \leq q$,
\begin{equation}
\nabla_j (g_{\vec{a}})(\vec{z})   = \rho_{\vec{a},t} \cdot \vec{e}_j + 2\pi a_j \cdot \nabla_j \whmu   \cdot \rho_{\vec{a},t} \cdot \vec{z}
\end{equation}
\begin{equation}
\| \nabla_j (g_{\vec{a}} - Id) \| \leq \|   \rho_{\vec{a},t} \cdot \vec{e}_j - \vec{e}_j \| + 2\pi \tn  \vec{a} \tn \cdot    B_1
\leq 2\pi (1 + B_1) \tn \vec{a} \tn
\end{equation}
 Thus, we may take  $C_1 = 2\pi (1 + B_1)$.
 The general case for $\|J | = r > 1$ proceeds similarly:
 \begin{eqnarray*}
\lefteqn{   \nabla_J (g_{\vec{a}} - Id)(\vec{z})      =    \nabla_J g_{\vec{a}} (\vec{z})      =   \nabla_{j_{r}} \circ \cdots \circ \nabla_{j_{1}} \left( g_{\vec{a}}(\vec{z}) \right) }  \\
& = & \nabla_{j_{r}} \circ \cdots \circ \nabla_{j_{2}} \left(\rho_{\vec{a},t} \cdot \vec{e}_{j_1} + 2\pi a_{j_1} \cdot \nabla_{j_1} \whmu   \cdot \rho_{\vec{a},t} \cdot \vec{z} \right) \\
& = & \nabla_{j_{r}} \circ \cdots \circ \nabla_{j_{3}} \left( 2\pi a_{j_2} \cdot \nabla_{j_2} \whmu   \cdot \rho_{\vec{a},t} \cdot \vec{e}_{j_1} + 2\pi a_{j_1} \cdot \nabla_{j_{2}}\nabla_{j_1} \whmu   \cdot \rho_{\vec{a},t} \cdot \vec{z}  + 2\pi a_{j_1} \cdot \nabla_{j_1} \whmu   \cdot  \nabla_{j_{2}}( \rho_{\vec{a},t} \cdot \vec{z}) \right)
\end{eqnarray*}
Observe that each term in parentheses in the last expression contains a factor of a component of $\vec{a}$, and the other factor involves  derivatives of $\whmu$.
 Thus, continuing on with the expansion, we obtain an expression where every term contains at least one factor $a_j$ for some $1 \leq j \leq q$, and the other factor involves derivatives of $\whmu$ times products of the components of  $\vec{a}$. By assumption, $\tn  \vec{a} \tn \leq 1$ so all products of components  of  $\vec{a}$ are likewise bounded above by $1$.
It  follows that there exists $C_{r}$  depending only on the estimates $B_{k}$ on the derivatives of $\whmu$ so that
 $\|  \nabla_J (g_{\vec{a}} - Id) \|' \leq C_{r}  \cdot \tn \vec{a}\tn$.
\endproof

\section{$C^r$-estimates}\label{sec-estimates}

The foliations $\whF_{\ell}$ constructed on $\whN_0 = \mT^k \times \mD^q$ in  section~\ref{sec-plug} are transverse to the factors $[\vec{x}] \times \mD^q$ for $[\vec{x}] \in \mT^k$, and so  can be alternately described in terms of their global holonomy maps, which define
  group actions $\whh_{\ell} \colon \mZ^k \to {\rm Diff}^{\infty}(\mD^q)$. We give explicit formulae for the maps $\whh_{\ell, j} = \whh_{\ell}(\vec{e}_j)$, $1 \leq j \leq k$, which are the generators of the $\mZ^k$-action. This yields criteria on the maps which are sufficient to guarantee that for each $j$,
the limit $\ds \whh_j = \lim_{\ell \to \infty} \whh_{\ell, j}$ is a $C^r$-diffeomorphism of $\mD^q$.

In the inductive construction of $\whF_{\ell +1}$, the   modification of  $\whF_{\ell}$ to  obtain $\whF_{\ell +1}$  is supported on  $\whN_{\ell}$, so the      holonomy maps of  $\whF_{\ell}$ and  $\whF_{\ell +1}$ agree  on $\whmK_{0} - \whmK_{\ell}$.
We obtain      formulae for  the holonomy of $\whF_{\ell +1}$ induced on $\whmK_{\ell}$ using that  the foliation $\whF_{\ell +1}$ is defined on the manifold  $N_{\ell}$ and then mapped to  $\whN_{\ell} \subset \whN_0 = \mT^k \times \mD^q$ via the map $\Psi_{\ell}$ defined by (\ref{eq-Psiell}). The map $\Psi_{\ell}$ ``twists" the product foliation on $N_{\ell}$ so that it agrees with the restriction of $\whF_{\ell}$ to $\whN_{\ell}$.
Estimates of the derivatives of differences $\whh_{\ell+1,j} - \whh_{\ell,j}$  on $\whmK_{\ell}$ then follow from   calculus.

The proof of the following key result is almost ``intuitively obvious'', as it is based on the effect of rescaling on the $C^r$-norm, but the dependence on the quantities $\tn \alpha_{\ell + 1}'(\vec{e}_j) \tn$ is perhaps less obvious.
\begin{prop}\label{prop-genest}
For all $r \geq 0$ and  $\ell \geq 0$,  then for   $\tn  \vec{a} \tn \leq 1$,
\begin{equation}\label{eq-estp}
\| \whh_{\ell +1 , j}  - \whh_{\ell , j} \|_r'  ~ \leq ~ C_{r} \cdot \e_{\ell}^{1-r} \cdot  \tn \alpha_{\ell + 1}'(\vec{e}_j) \tn
\end{equation}
where $C_r$ is defined by Lemma~\ref{lem-muest}, and  $\alpha_{\ell +1}' = \alpha_{\ell +1} \circ \Upsilon_{\ell}$ is defined below.
\end{prop}
\proof
Identify $\mD^q$ with the section $\whmK_0 = \Phi_0([\vec{0}] \times \mD^q_{\e_0}) = \Phi_0([\vec{0}] \times \mD^q)$.

Recall  that   the holonomy generators of $\F_{\ell +1}$ were
specified on $\pi_1(N_{\ell}, x_{\ell}) \cong \mZ^k$.

The inclusion  $\Psi_{\ell} \colon N_{\ell} \to \whN_{\ell} \subset \whN_0$ defined by  (\ref{eq-Psiell}) induces the
map as $\Phi_{\ell} $ defined by  (\ref{eq-Phi}),
$$\Phi_{\ell} \colon \mZ^k \cong    \G_{\ell} =  \pi_1(\whN_{\ell}, \whx_0) \subset \G_0 =  \pi_1(\whN_0, \whx_0) \cong \mZ^k$$
We  extend the holonomy of $\F_{\ell +1}$ from maps on the subgroup  $\G_{\ell} \subset \G_0$ to all of $\G_0$

First, extend
$\Phi_{\ell}$  to an isomorphism $\Phi_{\ell} \colon \mQ^k \to
\mQ^k$, which has inverse $\Upsilon_{\ell} \equiv \Phi_{\ell}^{-1}
\colon \mQ^k \to \mQ^k$ .
The isomorphism $\Upsilon_{\ell}$ is
represented by a matrix with rational entries, with l.c.d. $q_{\ell} = \det (\Phi_{\ell})$.

Intuitively, for $\gamma \in \pi_1(\whN_0, \whx_0)$ the rational number $\Upsilon_{\ell}(\gamma)$ is the lift of a path representing $\gamma$   to a ``fractional part'' of a closed path in the embedded torus
$\mT^k \cong  \whL_{\ell} \subset  \whN_{\ell} \subset  \whN_0 = \mT^k \times \mD_1^q$.
 The lift need not be a closed curve, unless
$\Upsilon_{\ell}(\gamma) \in \mZ^k$, but has initial and terminal points in the section $\mD^q$.

For each $\ell \geq 0$, $\alpha_{\ell +1} \colon \mZ^k \to \mQ^m$ admits a unique extension
$\alpha_{\ell +1} \colon \mQ^k \to \mQ^m$. Set   $\alpha_{\ell +1}' = \alpha_{\ell +1} \circ \Upsilon_{\ell}$.
The rational number  $\alpha_{\ell +1}'(\gamma)$ is the ``fractional rotation'' contributed by the
representation $\alpha_{\ell +1}$ along  the ``fractional loop'' $ \Upsilon_{\ell}(\gamma)$.
Define
\begin{equation}\label{eq-beta}
\beta_{\ell +1} \colon \mZ^k \to \mQ^m ~ , \quad \beta_{\ell +1} = \alpha_1' + \alpha_2' + \cdots + \alpha_{\ell +1}' ~,
\quad \beta_{\ell +1, j} = \beta_{\ell +1}(\vec{e}_j) \in \mQ^m
\end{equation}
\begin{equation}\label{eq-upsilon1}
\Theta_{\ell +1} \colon \mZ^k \to \SOq ~, ~ \Theta_{\ell +1}(\gamma) = \eb(\beta_{\ell +1}(\gamma)) ~  \text{for} ~ \gamma \in \mZ^k ~,  ~
\Theta_{\ell +1, j} = \Theta_{\ell +1}(\vec{e}_j) = \eb(\beta_{\ell +1, j})
\end{equation}

The next step is to give a formula for the centers $\{ \whz_{\ell +1, \gamma} \mid \gamma \in \mZ^k\}$ of the  disks of radius $\e_{\ell +1}$ which comprise $\whmK_{\ell +1} \subset \whmK_0$. We proceed by induction.

First note that for all $\vec{z} \in \mD^q_{\e_0'}  = \mD^q_{2/3}$ the holonomy action $h_1(\gamma)$ defined by (\ref{eq-holo1}) is equal to scalar multiplication by $ \Theta_1(\gamma)$. In particular, given the choice of generic vector $\vec{z}_1$  the centers of the $\e_1$-disks comprising $\whmK_1$ are the points
$$
\whz_{1,\gamma} =    \Theta_1(\gamma) \cdot \vec{z}_{1} \in \mD^q
$$
and the action of $\whh_1(\gamma')$ on $\whmK_1$ is also given by multiplication by $ \Theta_1(\gamma')$.

The centers of the disks  comprising $\whmK_2$ are obtained from the centers for $\whmK_1$ by adding on the rotation of the generic point $\vec{z}_2$ by the action of $\rho_2$.   In terms of the fundamental group $ \pi_1(\whN_0, \whx_0)$, these points are given by
\begin{equation}\label{eq-extend2}
\vec{z}_{2,\gamma} = \rho_{2}(\Upsilon_{1}(\gamma)) (\vec{z}_{2}) = \eb(\alpha_{2}'(\gamma)) \cdot \vec{z}_{2} \in \mD^q_{\e_{1}} ~ , ~ \text{for}~  \gamma \in \G_{1} \subset \mZ^k
\end{equation}
We use formula (\ref{eq-extend2}) to extend this action to all $\gamma \in \mZ^k$, but note that we must also multiply by  the term
 $\Theta_1(\gamma)$ which accounts for the rotation of $\whmK_1$ by $ \Theta_1(\gamma)$.
 Note that   $\Theta_1(\gamma)  = Id$ for $\gamma \in \G_1$ by definition, so this term does not appear in the
 formula (\ref{eq-extend2}). We then have, for all $\gamma \in \mZ^k$,
\begin{equation}
\whz_{2,\gamma} =  \whz_{1,\gamma} +  \eb(\alpha_{2}'(\gamma)) \eb(\alpha_{1}'(\gamma)) \cdot \vec{z}_{2} =
 \Theta_1(\gamma) \cdot \vec{z}_{1} +  \Theta_2(\gamma) \cdot \vec{z}_{2} \in \mD^q
\end{equation}
The action of $\whh_2(\gamma')$ on $\whmK_2$ is a linear isometry, where for $\whz_{2,\gamma} + \vec{v} \in \whmK_2$ with $\| \vec{v} \| \leq \e_2$ we have
\begin{equation}\label{eq-extend2c}
\whh_2(\gamma')(\whz_{2,\gamma} + \vec{v}) = \Theta_1(\gamma + \gamma') \cdot \vec{z}_{1} +  \Theta_2(\gamma + \gamma') \cdot (\vec{z}_{2} +\vec{v})
\end{equation}

For $\ell \geq 2$, with the choice of the generic points  $\vec{z}_{j+1} \in \mS^{q}_{\e_{j}/2} \subset \mB^{q}_{\e_{j}'} \cap \cO_{j +1}$ for the action $\rho_{j +1}$, for $1 \leq j < \ell$, we then have,  for all $\gamma \in \mZ^k$,
\begin{equation}
\whz_{\ell,\gamma} =   \Theta_1(\gamma) \cdot \vec{z}_{1} +  \Theta_2(\gamma) \cdot \vec{z}_{2} + \cdots +  \Theta_{\ell}(\gamma) \cdot \vec{z}_{\ell }  \in \mD^q
\end{equation}
Note that for $\gamma' \in \G_{\ell}$ we have   $\vec{z}_{\ell , \gamma + \gamma'} = \vec{z}_{\ell , \gamma}$.
In particular,
$\vec{z}_{\ell , \gamma'} = \vec{z}_{\ell , 0} = \vec{z}_{\ell}$.
It follows that
\begin{equation}\label{eq-Kell}
\whmK_{\ell} ~ = ~ \bigcup_{\gamma \in \mZ^k} ~ \mD^q_{\e_{\ell}}(\whz_{\ell , \gamma}) ~ \subset ~ \whmK_{\ell} \subset \whmK_0 = \mD^q
\end{equation}
The action of $\whh_{\ell}(\gamma')$ on $\whmK_{\ell}$ is a linear isometry, where for $\whz_{\ell,\gamma} + \vec{v} \in \whmK_{\ell}$ with $\| \vec{v} \| \leq \e_{\ell}$ we have
\begin{equation}\label{eq-extend2cc}
\whh_{\ell}(\gamma')(\whz_{\ell , \gamma} + \vec{v}) = \Theta_1(\gamma+ \gamma') \cdot \vec{z}_{1} + \cdots +  \Theta_{\ell}(\gamma+ \gamma') \cdot (\vec{z}_{\ell} + \vec{v})
\end{equation}

The last ingredient needed for the description of   the holonomy maps $\whh_{\ell +1,j}$ in terms of  $\whh_{\ell,j}$ is a family of
 rescaling maps for each of the disks in (\ref{eq-Kell}).  For each $\gamma \in \mZ^k$, set
\begin{equation}\label{eq-deflambda}
\lambda_{\ell, \gamma} \colon \mD^q \to \mD^q_{\e_{\ell}}(\whz_{\ell, \gamma}) ~ , ~ \lambda_{\ell, \gamma}(\vec{z}) =  \e_{\ell} \cdot \vec{z} + \whz_{\ell, \gamma}
\end{equation}

\medskip

We now give      formulae for the holonomy maps $\whh_{\ell}$.
First,   $\whh_{0,j} \colon \mD^q \to \mD^q$ is the identity map  for all $1 \leq j \leq k$,  as $\whF_0$ is the product foliation.

 Recall that $\whF_{\ell}$ and $\whF_{\ell + 1}$ are equal on $\whN_0 - \whN_{\ell}$, so $\whh_{\ell, j}$ and $\whh_{\ell +1 , j}$ agree on the set $\mD^q - \whmK_{\ell}$.  On the set $\whN_{\ell}$  recall that  the map $\Psi_{\ell} \colon N_{\ell} \to \whN_{\ell}$ defined in (\ref{eq-Psiell}) maps the product foliation on $N_{\ell}$ to the restriction of $\whF_{\ell}$ to $\whN_{\ell}$.
The foliation $\whF_{\ell+1}$ is obtained from $\whF_{\ell}$ by  ``twisting'' the product foliation on $N_{\ell}$ to a new foliation $\F_{\ell +1}$, and pushing the new foliation forward by $\Psi_{\ell}$.

For $\mD^q_{\e_{\ell}}(\whz_{\ell, \gamma}) \subset K_{\ell}$ and $\gamma' \in \mZ^k$, set
\begin{equation}
\whg_{\ell, \gamma}(\gamma')  =  \lambda_{\ell, \gamma} \circ g_{\alpha_{\ell + 1}'(\gamma')} \circ \lambda_{\ell, \gamma}^{-1} \colon \mD^q_{\e_{\ell}}(\whz_{\ell, \gamma})\to \mD^q_{\e_{\ell}}(\whz_{\ell, \gamma})
\end{equation}

Thus,  for   $\vec{v} \in \mD^q_{\e_{\ell}}(\whz_{\ell, \gamma}) \subset K_{\ell}$ and $\gamma' \in \mZ^k$, the   formula for the holonomy of $\whF_{\ell +1}$ is given by
\begin{equation} \label{eq-affine}
\whh_{\ell+1}(\gamma)(\vec{v})    =   \whh_{\ell, \gamma'} \circ   \whg_{\ell, \gamma}(\gamma')(\vec{v})
\end{equation}

Finally, we prove the estimate (\ref{eq-estp}). Restrict to a set  $\mD^q_{\e_{\ell}}(\whz_{\ell, \gamma}) \subset K_{\ell}$ then we have
\begin{equation}
\| \whh_{\ell +1 , j}  - \whh_{\ell , j} \|_r'   =      \|  \whh_{\ell, \vec{e}_j} \circ   \{ \whg_{\ell, \gamma}(\vec{e}_j) - Id\}  \|_r'
  =     \|  \whg_{\ell, \gamma}(\vec{e}_j) - Id  \|_r'
\end{equation}
 where we use that $\whh_{\ell, \vec{e}_j} $ acts via linear isometries on $\whmK_{\ell}$ so preserves $C^r$-norm. Then we have
  \begin{equation}
  \|  \whg_{\ell, \gamma}(\vec{e}_j) - Id  \|_r'
  =    \| \lambda_{\ell, \gamma} \circ  ( g_{\alpha_{\ell + 1}'(\vec{e}_j)} - Id ) \circ   \lambda_{\ell, \gamma}^{-1}  \|_r'
\end{equation}
The map $\lambda_{\ell, \gamma}$ is affine, so by the Chain Rule and Lemma~\ref{lem-muest} we obtain
   \begin{equation}
  \| \whh_{\ell +1 , j}  - \whh_{\ell , j} \|_r'
  =  \e_{\ell}^{1-r}  \cdot \|  ( g_{\alpha_{\ell + 1}'(\vec{e}_j)} - Id )    \|_r'
  \leq    \e_{\ell}^{1-r} \cdot C_{r} \cdot  \tn \alpha_{\ell + 1}'(\vec{e}_j) \tn
\end{equation}
This completes the proof of Proposition~\ref{prop-genest}.
\endproof

Proposition~\ref{prop-genest} gives $C^r$-norm estimates on the holonomy maps $h_{\ell}$ associated to a sequence of representations   $\whalpha = \{\alpha_{\ell} \colon \mZ^k \to \mQ^m \mid \ell \geq 0\}$  in terms of the quantities $\tn \alpha_{\ell +1}' \tn$  for each extended  representation   $\alpha_{\ell +1}' = \alpha_{\ell +1} \circ \Upsilon_{\ell}$.
If these norm estimates converge sufficiently rapidly, so that the sequence of maps $\{h_{\ell} \mid \ell \geq 0\}$ are Cauchy in the $C^r$-norm,  then their limit defines a $C^r$-action, and hence the foliations $\whF_{\ell}$ converge to a $C^r$-foliation $\whF$.

We introduce a quantity which measures   this ``total twisting'' for the   data $\whalpha' \equiv \{\alpha_{\ell}' \mid \ell \geq 0\}$, which depends upon both the    sequences of representations $\whalpha$ and associated bases $\{e^{\ell}_1, \ldots , e^{\ell}_k \}$  for the subgroups $\G_{\ell}$ they inductively  determine. Set
\begin{equation}
\tn \whalpha' \tn  = \sum_{\ell  = 0}^{\infty} ~  \tn \alpha_{\ell}'  \tn
\end{equation}

 There are three cases to consider, each with a distinct flavor: when $r = 0$, $r =1$ and $r = \infty$.

\begin{thm}\label{thm-C0}
Let  $\whalpha = \{ \alpha_{\ell} \colon \mZ^k \to \mQ^m \mid \ell \geq 1\}$ be an arbitrary  sequence of   representations.
Then there exists a sequence $\{\e_{\ell} > 0 \mid \ell \geq 0\}$ such that   the sequence of foliations $\{ \whF_{\ell} \mid \ell \geq 0\}$  converge in the $C^0$-topology to a foliation $\whF$ on $\mT^k \times \mD^q$ with minimal set $\cS_{\whalpha}$.
\end{thm}
\proof
By the   case $r =0$ of the estimate (\ref{eq-estp}) in Proposition~\ref{prop-genest}, we have
$$\| \whh_{\ell +1 , j}  - \whh_{\ell , j} \|_0  ~ \leq ~ 2\pi\cdot \e_{\ell} \cdot  \tn \alpha_{\ell + 1}'(\vec{e}_j) \tn$$
where   $C_0 = 2\pi$ by Lemma~\ref{lem-muest}. For each $\ell \geq 0$ choose
$$\e_{\ell +1} \leq \min \left\{  1/6^{\ell +1} ~ , 1/ (2^{\ell +1}    \tn\alpha_{\ell + 1}' \tn ) \right\}$$
Then the   holonomy maps $ \whh_{\ell , j}$   form a Cauchy sequence in the $C^0$-norm for each $1 \leq j \leq k$.
\endproof
This result is confirming the intuitively clear fact that   any sequence of coverings $p_{\ell} \colon L_{\ell} \cong \mT^k \to \mT^k$   can be realized as the minimal set of a $C^0$-foliation. This is in accord with the results of \cite{ClarkFokkink2004}.

Next, consider the case of $C^1$-embeddings:

\begin{thm}\label{thm-C1}
Let  $\whalpha = \{ \alpha_{\ell} \colon \mZ^k \to \mQ^m \mid \ell \geq 1\}$ be a  sequence of   representations such that    $\tn  \alpha_{\ell} \tn \leq 1$ for $\ell$ sufficiently large, and that bases $\{e^{\ell}_1, \ldots , e^{\ell}_k \}$  for the associated subgroups $\G_{\ell}$ have been chosen so that
\begin{equation}\label{eq-C1bounds}
\tn \whalpha' \tn  = \sum_{\ell  = 1}^{\infty} ~  \tn \alpha_{\ell}'  \tn  < \infty
\end{equation}
Then the sequence   $\{ \whF_{\ell} \mid \ell \geq 0\}$  converges in the $C^1$-norm  to a $C^1$-foliation $\whF$ on $\mT^k \times \mD^q$ with minimal set $\cS_{\whalpha}$.
\end{thm}
\proof
By Theorem~\ref{thm-C0} there exists a choice of diameters $\ve_{\ell} \to 0$
 such that  $ \whF_{\ell} \to \whF$ in the $C^0$-topology.
 Then, by the   case $r =1$ of the estimate (\ref{eq-estp}) in Proposition~\ref{prop-genest}, for $\ell$ sufficiently large, we have the estimate on the holonomy maps
$$\| \whh_{\ell +1 , j}  - \whh_{\ell , j} \|_1  ~ \leq ~  2\pi (1 + B_1) \cdot  \tn \alpha_{\ell + 1}'(\vec{e}_j) \tn$$
where    $C_1 = 2\pi (1 + B_1)$ by Lemma~\ref{lem-muest}, and $B_1 ~ = ~ \| \whmu \, \|_1'$. Then by  (\ref{eq-C1bounds}) the sequence $\{ \whh_{\ell , j} \mid \ell \geq 0\}$ is Cauchy in the $C^1$-norm for each $1 \leq j \leq k$.
\endproof
Again, the hypothesis (\ref{eq-C1bounds}) is  intuitively correct, based on the principle that rescaling a map does not change its $C^1$-norm. Thus, for the sequence $\{ \whF_{\ell} \mid \ell \geq 0\}$ to be Cauchy in the $C^1$-norm, it is required that the slopes of the modifications should be summable, which is  (\ref{eq-C1bounds}).

\medskip

For the smooth case,   we require that the foliations $\whF_{\ell}$ converge in the $C^r$-norm for all $r \geq 0$.  Naive intuition suggests that if the slopes of each successive modification tend to zero ``very quickly'' then the presentation $\cP = \{p_{\ell} \colon \mT^k \to \mT^k \mid \ell \geq 0\}$ can be realized as the compact leaves of a smooth foliation, limiting to an embedded solenoid homeomorphic to $\cS_{\cP}$.

On the other hand,  the estimate (\ref{eq-estp}) in Proposition~\ref{prop-genest}    for the $C^r$-norms $\| \whh_{\ell +1 , j}  - \whh_{\ell , j} \|_r'$   has  a factor   $\e_{\ell}^{1-r}$, and for $r > 1$ this tends to infinity as $\e_{\ell} \to 0$. Thus, a convergence criteria for $r > 1$ is more subtle than in the previous two cases $r=0$ and $r =1$. To the authors' knowledge, the results in the literature only discuss the   case $r \leq 1$ in detail, but do not address the more delicate issues of convergence for $r > 1$, even in the case of dimension-one solenoids.
We formulate below an existence theorem for a $C^{\infty}$ embedding, sufficient for our application to the proof of Theorem~\ref{thm-main-torus}.

\begin{thm}\label{thm-Cinf}
Let  $\whalpha = \{ \alpha_{\ell} \colon \mZ^k \to \mQ^m \mid \ell \geq 1\}$ be a  sequence of   representations such that    $\tn  \alpha_{\ell} \tn \leq 1$ for $\ell$ sufficiently large, and that radii $\{\e_1, \ldots, \e_{\ell}\}$  and bases $\{e^{\ell'}_1, \ldots , e^{\ell'}_k \}$  for the associated subgroups $\G_{\ell'}$ have been inductively chosen for $\ell' \leq \ell$.
Then assume that
\begin{equation}\label{eq-Cinfbounds}
  \tn \alpha_{\ell  +1}' \tn   \leq  (\e_{\ell})^{\ell}/(2^{\ell} \cdot \whC_{\ell +1})
\end{equation}
where $\whC_{\ell}$ is  defined by (\ref{eq-normmax}).
 Then for all $r \geq 0$, the family $\{\whF_{\ell} \mid \ell \geq 0\}$ converges in the  $C^{r}$-norm to a foliation $\whF$ on $\mT^k \times \mD^q$. Thus, $\whF$ is a smooth foliation with    minimal set $\cS_{\whalpha}$ which is a solenoid with presentation  $\cP$, formed by the coverings of $\mT^k$ associated to the tower of subgroups $\{\G_{\ell} \mid \ell \geq 0\}$.
\end{thm}
\proof
For each $r \geq 0$, estimate (\ref{eq-estp}) and hypothesis (\ref{eq-Cinfbounds}) imply that
$\{ \| \whh_{\ell +1 , j}  - \whh_{\ell , j} \|_r \mid \ell \geq 0 \}$ is a Cauchy sequence. Thus, the limit foliation $\whF$ is $C^r$.
\endproof

We make two remarks concerning the hypotheses (\ref{eq-Cinfbounds}). First, it implies that the bounds $ \tn \alpha_{\ell  +1}' \tn$ on the slopes of each modification of $\whF_{\ell}$ to obtain $\whF_{\ell +1}$ tend to zero \emph{very rapidly}, at least exponentially.

Secondly, the choice of each subsequent radius $\e_{\ell +1}$ is made after the choice of $\alpha_{\ell}$, while the bound (\ref{eq-Cinfbounds}) for $\ell +1$ depends only on the prior choices of radii, $\{\e_1, \ldots, \e_{\ell}\}$. Thus, one can use this latitude to successively choose the representations $\alpha_{\ell}$ so that  (\ref{eq-Cinfbounds}) is satisfied for all $\ell \geq 0$.  These ideas are the basis for our constructions in section~\ref{sec-foliations}.

\section{$C^r$-embedded $1$-dimensional solenoids}\label{sec-flows}

We consider first   the traditional case of dimension-one solenoids, where many of the   estimates in section~\ref{sec-estimates} are greatly simplified.
This is the most familiar and intuitive case, for which there is an extensive literature (see, for example \cite{BowenFranks1976, Gambaudo1991, GST1994, GT1990, Kan1986, MM1980, Moreira2004, EThomas1973}).
Of course, the main point of this paper is to give an explicit  construction which yields   smoothly embedded  higher dimensional solenoids, but examining this simplest case first  illustrates the   steps of the induction in sections~\ref{sec-plug} and \ref{sec-estimates}. The higher dimensional cases are considered subsequently.

Let $k=1$ and $q=2$. The case for $q > 2$ is handled in a similar fashion, and is left to the reader.

Assume there is given a presentation $\cP = \{ p_{\ell} \colon \mS^1 \to \mS^1 \mid \ell \geq 1\}$ as in Definition~\ref{def-solenoid}. We assume the maps are oriented, thus $p_{\ell}$ is determined up to isotopy by its degree $m_{\ell} > 1$. Thus, we will assume the maps are in standard form, with $p_{\ell}(z) = z^{m_{\ell}}$, where $z$ is the complex coordinate for $\mS^1 \subset \mC$. The inverse limit solenoid $\cS_{\cP}$ is then determined up to homeomorphism by the set $\{m_1, m_2, \ldots\}$.

For each $\ell > 0$ we require a   representation $\rho_{\ell} \colon \mZ \to \SO(2)$ with kernel $\G_{\ell} = m_{\ell} \cdot \mZ$, and its lift to a   representation  $\alpha_{\ell} \colon \mZ \to \mQ$. This is equivalent to the choice of $a_{\ell} = \alpha_{\ell}(1)   \in \mQ$ such that $a_{\ell} \mod \mZ$ is a  root-of-unity of order   $m_{\ell}$ in $\mQ/\mZ$. For the following, we make the ``standard''  choice  of  $a_{\ell} = 1/m_{\ell}$ and so  $\alpha_{\ell}(i) = i/m_{\ell} \in \mQ$.
Note that for each $\ell > 0$,  there is an infinite number of possible choices at this stage of the construction. On the other hand, there is a unique  choice  $e_1^{\ell} = m_{\ell} \cdot 1$ of the oriented  generator of  $\G_{\ell} = m_{\ell} \cdot \mZ$ which is the kernel of $\rho_{\ell}$. 
Introduce the  notation   $ \alpha_{\ell}(r) = \alpha_{\ell}(e_1^{\ell} )(r)$.

We assume there is given an infinite string of pairs $\{(m_{\ell} , a_{\ell}) \mid \ell = 1,2, \ldots \}$, where the positive integers $\{m_{\ell} \}$ determine  the homeomorphism type of the solenoid $\mS_{\cP}$, and  the rational numbers $\{a_{\ell} \}$ determine the embedding of $\mS_{\cP}$ into $\mS^1 \times \mD^2$.

Following the notation of Section~\ref{sec-estimates}, the  map  $\phi_{\ell}  \colon \mR \to \mR$ is   given by $\phi_{\ell}(r) = m_{\ell} \,  r$,   inducing the  diffeomorphism $\pb_{\ell} \colon \mR/\mZ \to \mR/m_{\ell} \, \mZ$.
The map $\Phi_{\ell} \colon \mZ \to \mZ$ is then given by   $\Phi_{\ell}(r) =d_{\ell} \cdot r$, where $d_{\ell} =  m_1 \cdots m_{\ell}$.
The inverse
 $\Upsilon_{\ell} \colon \mQ \to \mQ$  is simply $\ds \Upsilon_{\ell}(s) = s/d_{\ell}$.
Then
\begin{equation}\label{eq-q=1ests}
\alpha_{\ell +1}'(r) = \alpha_{\ell +1} \circ \Upsilon_{\ell}(r) =  a_{\ell +1} \, r/ d_{\ell} =  r/d_{\ell+1}
\end{equation}
\begin{equation*}
\beta_{\ell} (r)    =    \left\{ \frac{1}{m_1} +
\frac{1}{m_1 m_2} + \cdots + \frac{1}{m_1 m_2 \cdots m_{\ell}} \right\}  \cdot r  = \beta_{\ell,1} \cdot r
\end{equation*}

The choices of the radii $\e_{\ell}$ are dictated by the conditions (\ref{eq-diskell}).
Given the choice of a   generic point $\vec{v}_{\ell} \in \mD^2_{\e_{\ell}}$ with $\|\vec{v}_{\ell}\| = \e_{\ell}/2$,
 the distance between the translates of $\vec{v}_{\ell}$ by the rotation group of order $m_{\ell} \geq 2$ is bounded below by
$\e_{\ell}/4m_{\ell}$.
Thus, it suffices to  require that
$\e_{\ell +1} \leq \e_{\ell}/8 m_{\ell}$ in order to ensure that the translates $\mD^2_{\e_{\ell +1}} \subset \mD^2_{\e'_{\ell}}$   under the rotation group of order $m_{\ell +1}$ are disjoint, and also to satisfy the condition that the action of $\rho^{{\ell +1}}_t$ is affine on the disks. Inductively, we see this implies
$\ds
\e_{\ell} \leq 1/8^{\ell} d_{\ell}
$ is a sufficient condition on the radii for the disks to be disjoint.

The hypotheses of  Theorem~\ref{thm-C0} are satisfied for all choices as above, including for all choices of   $a_{\ell} = \alpha_{\ell}(1)   \in \mQ$ such that $a_{\ell} \mod \mZ$ is a  root-of-unity of order   $m_{\ell}$ in $\mQ/\mZ$. Thus, we see that every ``standard'' $1$-dimensional solenoid  admits an embedding, along with its given presentation,  into a $C^0$-foliation $\whF$ of $\mS^1 \times \mD^2$.

Next, consider the embedding problem for a given presentation $\cP$ of a solenoid $\cS \cong \cS_{\cP}$ into a $C^1$-foliation. 
The hypotheses of Theorem~\ref{thm-C1} are simply that $a_{\ell} \leq 1$ for $\ell$ sufficiently large, and that
 \begin{equation}\label{eq-C1embedbds}
\tn \whalpha' \tn  = \sum_{\ell  = 1}^{\infty} ~  \tn \alpha_{\ell}'  \tn  ~ = ~ \sum_{\ell  = 1}^{\infty} ~ \frac{1}{m_1 \cdots m_{\ell}}   ~  < \infty
\end{equation}
It is assumed that each $m_{\ell} > 1$, so condition (\ref{eq-C1embedbds}) is satisfied by all choices as above. Thus, given the simple assumption that the  embedded tube $\whN_{\ell +1} \subset \whN_{\ell}$ has slope less than $1$ for $\ell$ sufficiently large, it follows that every   $1$-dimensional solenoid  $\cS_{\cP}$ admits an embedding, along with its given presentation $\cP$,  into a $C^1$-foliation $\whF$  of $\mS^1 \times \mD^2$.

Lastly, for the     case of a smooth embedding, the hypotheses of Theorem~\ref{thm-Cinf} are not automatically satisfied for every presentation of a given solenoid. Recall the hypotheses of the theorem are   sufficient conditions for the existence of an embedding, which are that  an  inductive condition be satisfied:
\begin{equation}\label{eq-inductest3a}
  \tn \alpha_{\ell  +1}' \tn   \leq  (\e_{\ell})^{\ell} / (2^{\ell} \cdot \whC_{\ell +1}\})
\end{equation}
where the constant   $\whC_{\ell +1}$ is  defined by (\ref{eq-normmax}). Choose $\e_{\ell} = 1/8^{\ell} d_{\ell}$, then this simplifies to the condition
\begin{equation}\label{eq-inductest3b}
m_{\ell +1}  \geq \sqrt{2^{\ell} 8^{\ell^2} (d_{\ell})^{\ell -1} \cdot \whC_{\ell +1}   }
\end{equation}

Given a presentation $\cP$ for a solenoid $\cS \cong \cS_{\cP}$, there exists a sub-presentation $\cP'$ of $\cP$ for which condition (\ref{eq-inductest3a})  is satisfied. Thus, we see that the homeomorphism class of each $1$-dimensional solenoid admits an embedding into a smooth foliation of codimension $2$. The differences between the $C^1$-embedding condition (\ref{eq-C1embedbds}) and the $C^{\infty}$-embedding condition (\ref{eq-inductest3b}), is a restriction on the lengths of the approximating periodic orbits for the flow in which the minimal set is realized as a minimal set.

  Note that  the orders $\{m_{\ell} \mid \ell \geq 1\}$ of the holonomy groups of the inserted disks tend to infinity quite rapidly. 
It is also worth noting that, none the less,  there are an uncountable number of sequences  $\{m_{\ell} \mid \ell \geq 1\}$  which satisfy the inductive criteria (\ref{eq-inductest3b}).

\section{$C^r$-embedded higher-dimensional solenoids}\label{sec-foliations}

Recall that the given a compact manifold $M_0$ of dimension $n$ and surjection $\pi_1(M_0, x_0) \to \mZ^k$, the construction of a $C^r$-foliation  of $\mT^k \times \mD^q$ with a solenoidal minimal set yields a   $C^r$-foliation of   $M_0 \times \mD^q$ with   a solenoidal minimal set. In this section, given a presentation $\cP$ for a minimal set $\cS_{\cP}$ over $\mT^k$, for $k \geq 2$, we  consider the problem of constructing $C^r$-foliations  of $\mT^k \times \mD^q$ with   minimal sets homeomorphic to $\cS_{\cP}$.

Let $\cS$ be a solenoid with    base space $\mT^k$ and presentation
 $\cP = \{ p_{\ell} \colon \mT^k \to \mT^k \mid \ell \geq 0\}$.
By the results of McCord, the solenoid is determined up to homeomorphism by the Cantor group bundle structure. Thus, we can assume that each map  $p_{\ell}$ is affine  with $p_{\ell}([\vec{0}]) = [\vec{0}]$, and there are   inclusion maps
 $$p_{\ell *} \colon \mZ^k \cong \pi_1( \mT^k, [\vec{0}]) \to  \mZ^k \cong \pi_1( \mT^k, [\vec{0}])$$
 with each image   a proper subgroup of finite index.
 Define
 $\pi_{\ell} = p_1 \circ \cdots \circ  p_{\ell} \colon \mT^k \to \mT^k$
 and set
$\G_{\ell} = {\rm Image} \left\{\pi_{\ell *} \colon \mZ^k \to \mZ^k \right\}$. Then $\G_{\ell}$ is a proper subgroup of finite index, hence is torsion-free of rank $k$.   The presentation $\cP$ defines a descending chain
\begin{equation}\label{eq-chainL}
\G_{\cP} \equiv \left\{ \mZ^k = \G_0 \supset \G_1 \supset \G_2 \supset \cdots \supset  \G_{\infty} \right\}
\end{equation}
 where each inclusion $\G_{\ell + 1} \subset \G_{\ell}$ has finite index greater than one for all $\ell \geq 0$,
and $\ds \G_{\infty}  = \bigcap_{\ell \geq 0} ~ \G_{\ell}$.

Our approach to constructing an embedding for the solenoid $\cS_{\cP}$ into $\mT^k \times \mD^q$ via flat bundles as in section~\ref{sec-plug}, is to obtain representations of each of the quotient groups  $\fG_{\ell} \cong \G_{\ell}/\G_{\ell +1}$ into $\mT^m \subset \SO(q)$.
In  section~\ref{sec-plug},  $\fG_{\ell} =  \Delta_{\ell}/\mZ^m \subset \mQ^m/\mZ^m$ was defined as the image of the representation
 $\rho_{\ell} = \rho^{\alpha_{\ell}} \colon \mZ^k \to \mT^m$.
As shown is section~\ref{sec-estimates}, the smoothness of the resulting foliation $\whF$ on  $\mT^k \times \mD^q$ depends upon the properties of these representations, as well as the ``positioning'' of the subgroups $\G_{\ell} \subset \G_0 = \mZ^k$. The next part of our analysis of a tower $\G_{\cP}$ analyzes this process in more detail.

In addition to the above data, the construction of an  embedding of $\cS_{\cP}$ also  requires    choices of the representations  $\rho_{\ell} \in {\rm Hom}(\fG_{\ell} , \mT^m)$, for $\ell \geq 1$, and   lifts $\alpha_{\ell} \colon \G_{\ell} \to \mQ^m$ such that    the estimates of section~\ref{sec-estimates} are satisfied.
The   ``Simultaneous Bases Theorem'' for free abelian groups, as given for example    in  \cite[Theorem~10.21]{Rotman1995}, provides one method for choosing the data  $\alpha_{\ell}$ and $\rho_{\ell}$.

 \begin{thm}[Simultaneous Bases Theorem] \label{thm-SBT}
 Let $H$ be a subgroup of finite index in a free abelian group $F$ of finite rank $k$. Then there exist bases $\{f_1, \ldots , f_k\}$ of $F$ and $\{h_1, \ldots , h_k\}$ of $H$ such that $h_i = m_i \cdot f_i$ for all $1 \leq i \leq k$, where $m_i \geq 1$.
 \end{thm}

We   apply  Theorem~\ref{thm-SBT} to the chain $\G_{\cP}$.
For all $\ell \geq 0$, there exists a basis $\bB_{\ell} \equiv \{g^{\ell}_1, \ldots , g^{\ell}_{k}\} \subset  \G_{\ell}$
and positive  integers $\vec{m}_{\ell +1} = (m_{\ell+1, 1}, \ldots, m_{\ell+1, k})$ such that
 $\vec{m}_{\ell +1} \cdot \bB_{\ell} = \{  m_{\ell+1, i} \cdot g^{\ell}_i \mid 1 \leq i \leq k \}$ is a basis for $\G_{\ell +1}$.
% The product $d_{\ell +1} = m_{\ell+1, 1} \cdots m_{\ell+1, k_1} > 1$  is the index of the subgroup $\G_{\ell +1}$ in $\G_{\ell}$.

  Define a sequence of ``standard representations''  $\whalpha = \{\alpha_{\ell} \colon \mZ^k \to \mQ^k \}$ by
 $\alpha_{\ell}(\vec{e}_j) = \frac{1}{m_{\ell+1, j}} \cdot \vec{e}_j$. Then for $q = 2k$, we obtain the    associated representations $\rho_{\ell} = \rho^{\alpha_{\ell}} \colon \mZ^k \to \SO(q)$ as in (\ref{eq-defrho}). The basis  $\bB_{\ell}$ defines an isomorphism $\phi_{\ell} \colon \mZ^k \to \G_{\ell}$ and the composition $\rho_{\ell} \circ \phi_{\ell}^{-1} \colon \G_{\ell} \to \SO(q)$ has kernel $\G_{\ell +1}$.

The sequence of data $\{(\alpha_{\ell}  , \bB_{\ell}) \mid \ell \geq 0\}$ as above is called a \emph{standard representation} for the chain $\G_{\cP}$.  Note that this sequence of data need not be uniquely determined by the chain $\G_{\cP}$, as the  choice of simultaneous basis $\bB_{\ell}$ for each pair $\G_{\ell +1} \subset \G_{\ell}$ need not be unique.

We use this standard representation to show that every presentation $\cP$ embeds in a $C^0$-foliation of codimension $q = 2k$. The constructions are easily modified to embed in any codimension $q \geq 2k$: 
 \begin{prop}\label{prop-embedC0}
Let $\cP$ be a presentation of the solenoid $\cS_{\cP}$ over the base space $\mT^k$, and let $q = 2k$. Then there exists a $C^0$-foliation $\whF$ of $\mT^k \times \mD^q$ such that there is an embedding of   $\cP$ into $\whF$. In particular, $\cS_{\cP}$ is homeomorphic to a minimal set $\fM$ for $\whF$.
\end{prop}
 \proof
We assume a standard representation  $\{(\alpha_{\ell}  , \bB_{\ell}) \mid \ell \geq 0\}$ has been chosen for $\cP$.
Each  basis  $\bB_{\ell}$ defines an isomorphism $\phi_{\ell} \colon \mZ^k \to \G_{\ell}$ and the composition $\rho_{\ell} \circ \phi_{\ell}^{-1} \colon \G_{\ell} \to \SO(q)$ has kernel $\G_{\ell +1}$.

Then choose a sequence of radii $\e_{\ell}$ as in  the proof of Theorem~\ref{thm-C0}. We obtain a sequence of foliations $\{\whF_{\ell} \mid \ell \geq 0\}$ on $\mT^k \times \mD^q$   converging $C^0$ to a foliation $\whF$ such that $\cP$ embeds into $\whF$ with minimal set $\fM$ homeomorphic to $\cS_{\cP}$.
 \endproof

It is known that a non-trivial solenoid cannot be embedded as a codimension $1$ submanifold, due to cohomology considerations (see \cite{ClarkFokkink2004, Prajs1990}), so in particular cannot be embedded as a minimal set for a codimension $1$ foliations.  It is unknown if the conclusion of Proposition~\ref{prop-embedC0} can possibly be extended for some  $1 < q < 2k$.  More likely,  it may be possible to prove the non-embedding in this range,  possibly  based on properties of higher order linking invariants of an embedding, modeled on the classical linking invariants, discussed  as   in \cite{ClarkSullivan2004, Gambaudo2006}.

There is no a priori control of the size of the integers $\{m_{\ell+1, j} \mid \ell \geq 0\}$ which arise in a standard representation of $\cP$, so the estimates (\ref{eq-C1bounds}) and (\ref{eq-Cinfbounds}) need not be satisfied for a standard representation, and the embedding of Proposition~\ref{prop-embedC0} need not be $C^1$.
In general, to obtain a $C^r$-embedding  of the solenoid  $\cS_{\cP}$, it is necessary to consider sub-presentations of $\cP$, or equivalently a sequence $0 = \ell_0 < \ell_1 < \ell_2 < \cdots$ which determines the subchain
$$\G_{\cP'} = \{\mZ^k = \G_0  \supset \G_{\ell_1} \supset \G_{\ell_2} \supset \cdots \}$$
We formulate an algebraic  condition on a chain $\G_{\cP'}$ such that  $\cS_{\cP} \cong  \cS_{\cP'}$ admits a $C^r$-embedding.

For $\ell \geq 0$,
for each pair $\G_{\ell +1} \subset \G_{\ell}$ for  $\ell \geq 0$, let  $\bB_{\ell} \equiv \{g^{\ell}_1, \ldots , g^{\ell}_{k}\} \subset  \G_{\ell}$   be a  simultaneous basis  with associated degree vectors $\vec{m}_{\ell +1}$.
Define
 $S_{0}  \in \SL(\mZ^{k})$ to be the matrix such that $S_0 \cdot \vec{e}_j = g^{0}_1$ where   $\{\vec{e}_1, \ldots , \vec{e}_k\}$ is the standard basis for $\mZ^k$. Set    $m (0; i) =1$ for $1 \leq i \leq k$, so $\vec{m}_0 = (1,1, \ldots , 1) \in \mZ^k$.

For each $\ell \geq 0$, the group $\G_{\ell +1}$ has  bases, $\vec{m}_{\ell +1} \cdot \bB_{\ell}$ and $\bB_{\ell +1}$, and so there exists   $S_{\ell +1}  \in \SL(\mZ^{k})$   such that $S_{\ell +1}   \cdot  ( m_{\ell +1, i} \cdot  g^{\ell}_i) = g^{\ell +1}_i$. The descending chain of subgroups  $\G_{\cP}$ is then completely determined by the data
 \begin{enumerate}
\item a sequence of matrices  $S_{\ell} \in \SL(\mZ^{k})$ for $\ell \geq 0$
\item degree vectors $\vec{m}_{\ell} = (m_{\ell, 1}, \ldots, m_{\ell, k})$  with $m_{\ell, i} \geq 1$, for $\ell \geq 1$, $1 \leq i \leq k$
\end{enumerate}
Introduce the matrix $\whS_{\ell}$ which defines the \emph{total change of framing} (or the ``total twisting'') in passing from the simultaneous framings $\cB_0$ for $\G_0$ and $\cB_{\ell}$ for $\G_{\ell}$:
\begin{equation}\label{eq-totalbase}
\whS_{\ell} = S_{\ell}   \cdot  S_{\ell -1} \cdots S_{0}  \in \SL(\mZ^{k})
\end{equation}
Also, introduce the diagonal matrices with entries to product of multipliers:
\begin{equation}\label{eq-totalexp}
 \whM_{\ell} = {\rm Diag}(\whm (\ell; 1), \ldots , \whm (\ell; k)) ~ , ~   \whm (\ell; i) = m_{\ell, i}  \cdot m_{\ell -1, i}  \cdot m_{1, i}
\end{equation}
Then we have the fundamental identity:
\begin{lemma}\label{lem-baseschain}
For $\ell \geq 0$, $\bB_{\ell} = \whS_{\ell}  \cdot  \whM_{\ell}  \cdot \bB_0$.
\end{lemma}
\proof
We must show that $g^{\ell}_i = \whS_{\ell} \cdot  \whm (\ell; i) \cdot \vec{e}_i$ for each $1 \leq i \leq k$.
This follows by noting that  $S_{\ell +1}   \cdot  m_{\ell +1, i} \cdot  g^{\ell}_i  = g^{\ell +1}_i$ and applying induction.
\endproof

We adapt this data to the subchain $\G_{\cP'}$ defined by $\G_j' \equiv  \G_{\ell_j}$.  For each $j \geq 1$,  set
\begin{equation}\label{eq-submatrix1}
 S_{j}' \equiv  S_{\ell_j}   \cdot  S_{\ell_j - 1} \cdots S_{\ell_{j-1} +1} \in \SL(\mZ^k) ~  , \quad    m'_{j; i} = m_{\ell_j , i}  \cdot m_{\ell_j -1, i}  \cdot m_{\ell_{j-1} +1 , i}
\end{equation}
\begin{equation}\label{eq-submatrix2}
\whS_{j}' = \whS_{\ell_j} \in \SL(\mZ^k) ~ , \quad   \whM_j' = \whM_{\ell_j} ~ , \quad \bB'_j = \bB_{\ell_j}
\end{equation}
Then  Lemma~\ref{lem-baseschain} directly implies that $\bB_{j}' = \whS_{j}' \cdot  \whM_j' \cdot \bB_0$.

 Define a sequence of representations $\alpha_{j}' \colon \mZ^k \to \mQ^k$ by
 $\alpha_j'(\vec{e}_i) = \frac{1}{m_{j+1, i}'} \cdot \vec{e}_i$  for $1 \leq i \leq k$, with associated representations $\varrho_{j} = \rho^{\alpha_{j}'} \colon \mZ^k \to \SO(q)$ as in (\ref{eq-defrho}). The basis  $\bB_{j}$ defines an isomorphism $\phi_{j}'  \colon \mZ^k \to \G_{j}'$ and the composition $\varrho_{j} \circ (\phi_{j}')^{-1} \colon \G_{j}' \to \SO(q)$ has kernel $\G_{j +1}'$.

 Let $\whalpha'$ denote the data $\{( S_{j}' ,  \vec{m}'_{j} \mid j \geq 1\}$ which call a standard representation of $\G_{\cP'}$

 Recall that $ \Phi_{\ell} \colon \mZ^k \to \mZ^k$ with image $\G_{\ell}$ was defined in (\ref{eq-Phi}). In terms of the above, we have:
 \begin{lemma}\label{lem-subPhi}
For $\ell \geq 1$, $\Phi_{j}' = \Phi_{\ell_j}$ is represented by the matrix $ \whS_{j}' \cdot  \whM_j'$. \hfill $\Box$
\end{lemma}

Let $\Upsilon_j'$ denote the inverse to $\Phi_j'$ which is represented by the matrix
$(\whS_j')^{-1}  \in \SL(\mZ^k)$.

Set   $\alpha_{j +1}'' = \alpha_{j +1}' \circ \Upsilon_{j}' \colon \mZ^k \to \mQ^k$.

Also,  set  $\whm_*(j) = \min \left\{  \whm(\ell_{j} , i) \mid 1 \leq i \leq k \right\} $ and
 $\whm^*(j) = \max \left\{  \whm(\ell_{j} , i) \mid 1 \leq i \leq k \right\}$.

Note that $\alpha_{j +1}''(\vec{e}_i) = \whm(\ell_{j+1} , i)^{-1} \cdot (\whS_j')^{-1} \cdot \vec{e}_i$ so that
\begin{equation}\label{eq-q=kests}
\tn \alpha_{j+1}''  \tn   =   \max \left\{ \tn \alpha_{j +1}''(\vec{e}_i)  \tn  \mid 1 \leq i \leq k \right\}   \leq    \frac{1}{  \whm_*(j+1)} \cdot \| (\whS_j')^{-1} \|
\end{equation}

Note that the quantity $\tn \alpha_{j+1}''  \tn $  depends only on the algebraic structure of the sub-chain $\G_{\cP'}$ and the choice of the simultaneous bases at each stage. Note the  analogy between the estimate (\ref{eq-q=kests}) and that obtained from (\ref{eq-q=1ests})  for the case $k=1$. In the estimate  (\ref{eq-q=1ests}), the matrix $\whS_j' = \pm Id$ while the term $d_{\ell +1}$ is the diagonal entry of the $1 \times 1$ matrix $\whM_{\ell}$, which is $\whm(\ell_{\ell +1} , 1)$.

Now define a quantity which depends only only on the algebraic structure of the sub-chain $\G_{\cP'}$
\begin{equation}\label{eq-C1algbounds}
\tn \G_{\cP'} \tn  = \inf \left\{ \sum_{j  = 0}^{\infty} ~ \tn \alpha_{j+1}''  \tn  \mid \whalpha' ~ \text{is a standard presentation for} ~ \cP' \right\}
\end{equation}

The following now follows from our previous constructions and results. Note that the $C^1$-norm estimate (\ref{eq-C1bounds}) is independent of the radii $\e_{\ell_j}$ chosen, so they do not appear in the hypotheses.
  \begin{prop}\label{prop-embedC1}
Let $\cP'$ be a sub-presentation of the solenoid $\cS_{\cP}$ over the base space $\mT^k$, and let $q = 2k$. If $\tn \G_{\cP'} \tn < \infty$, then   there exists a $C^1$-foliation $\whF$ of $\mT^k \times \mD^q$ such that there is an embedding of   $\cP'$ into $\whF$. In particular, $\cS_{\cP}$ is homeomorphic to a minimal set $\fM$ for $\whF$.
\end{prop}

 Likewise, the above estimates for a sub-presentation $\cP'$ can be applied to obtain smooth embeddings.  In this case, it is necessary to also have an estimate for the radii $\e_j$ appearing in the  estimates (\ref{eq-Cinfbounds}). To obtain these, we note that at each stage of the construction of the plug $\whN_j$ in section~\ref{sec-plug}, the representation $\rho_{j+1}'$ is a product of $1$-dimensional representations. Thus, the \emph{minimum} distance between points in an orbit  of $\rho_{j+1}'$ can be estimated as in section~\ref{sec-flows}, to obtain that it suffices to let
  $\e_{j+1} = 1/8^j \whm^*(j)$, where now we use the maximum of the exponents, $\whm^*(j)$, for our estimate.

The estimate (\ref{eq-Cinfbounds}) assumes that
$\ds   \tn \alpha_{j  +1}'' \tn   \leq  (\e_{j})^{j}/(2^{j} \cdot \whC_{\ell_j +1})$ so by (\ref{eq-q=kests}) it suffices to require    that
$$ \frac{1}{  \whm_*(j+1)} \cdot \| (\whS_j')^{-1} \|  \leq   (\e_{j})^{j}/(2^{j} \cdot \whC_{\ell_j +1})$$
Then the following now follows from our previous constructions and results:
  \begin{prop}\label{prop-embedCinf}
Let $\cP$ be a presentation of the solenoid $\cS$ over the base space $\mT^k$, let $q = 2k$. Suppose that there exists a sub-presentation $\cP'$ of $\cP$ and standard representation $\whalpha'$ for $\cP'$ such that
\begin{equation}\label{eq-Calginfbounds}
\whm_*(j+1)   \geq  \frac{(8^j \whm^*(j))^{j} \cdot 2^{j} \cdot \whC_{j +1}}{ \whm_*(j+1)  \| (\whS_j')^{-1} \| }
\end{equation}
holds for all $j$ sufficiently large.
Then   there exists a $C^r$-foliation $\whF$ of $\mT^k \times \mD^q$ such that there is an embedding of   $\cP'$ into $\whF$. In particular, $\cS_{\cP}$ is homeomorphic to a minimal set $\fM$ for $\whF$.
\end{prop}

The unwieldy formula (\ref{eq-Calginfbounds}) can be viewed as a prescription for constructing towers $\cP$ for which the corresponding inverse limit $\cS_{\cP}$ embeds into a smooth foliation of $\mT^k \times \mD^q$. Start with the standard group $\G_0 = \mZ^k$, and choose a basis $\bB_0$ for $\G_0$. Then choose any vector $\vec{m}_1$ such that   all entries   $m_{1,i} \geq 1$. Then let $\G_1$ be the subgroup generated by $\vec{m}_1 \cdot \bB_0$. Choose a basis $\bB_1$ for $\G_1$, or equivalently a matrix $S_1 \in \SL(\mZ^k)$ so that
$\bB_1 = S_1 \cdot \vec{m}_1 \cdot \bB_0$. Then repeat as long as desired, though at some stage in the construction, it is necessary to require that subsequent choices of the multipliers $\vec{m}_{\ell}$ are sufficiently large so that (\ref{eq-Calginfbounds}) holds for all $j \to \infty$. The resulting foliation $\whF$ will then be smooth, by Theorem~\ref{thm-Cinf}.

Note that if all choices of the multiplier vectors $\vec{m}_{\ell}$ in the above construction scheme have all components sufficiently large, then the resulting foliation $\whF$ will be arbitrarily close to the product foliation on $\mT^k \times \mD^q$.

The claims of Theorem~\ref{thm-main-fol}  now follow from these remarks and our previous constructions.

Finally, we indicate how Corollary~\ref{cor-uncountable} is derived from Theorem~\ref{thm-main-fol}.
    First,  observe that for an $n$-dimensional solenoid $\cS$ as presented in
Definition~\ref{def-solenoid}, the $n$-th \v{C}ech cohomology group
$\check{H}^n(\cS ;\mZ)$ is isomorphic to the direct limit
$\underrightarrow{\lim}~\{d_{\ell} \colon \mZ \rightarrow \mZ\}$,
where $d_{\ell}$ is the degree of the covering map $p_{\ell}$. This
is a torsion-free group of rank one. Without loss of generality, one
may assume that any such group is presented as
$\underrightarrow{\lim}~\{p_{\ell} \colon \mZ \rightarrow \mZ\}$,
$\ell \geq 1$, where each map $p_{\ell}$ is multiplication by a
prime $p_{\ell}$.
 According to Baer's
classification of such groups \cite{Baer1937, Kechris2000}, two such
groups, which are determined by the sequences of primes $P=(p_1,p_2,\dots)$ and
$Q=(q_1,q_2,\dots)$, are isomorphic if and only if it is possible to
remove finitely many of the terms from each of the sequences $P$ and
$Q$ to obtain new sequences $P'$ and $Q'$ in such a way that each
prime number occurs with the same cardinality in $P'$ and $Q'$, see
e.g. McCord~\cite{McCord1965}.
Thus, by choosing these degrees
appropriately, and tending to infinity very rapidly,
it is clear that there is an uncountable number of
topologically distinct solenoids based on a given $L_0$ which embed smoothly as minimal sets of $C^{\infty}$-foliations,
as in Theorem~\ref{thm-main-fol}.

At each stage of our construction, the foliated plug constructed in section~\ref{sec-plug} has a central disk which is invariant, and on which the leaves remain compact without holonomy under all subsequent modifications. Thus, the limiting foliation $\cH_1 = \whF$ contains open saturated sets, foliated by compact leaves with trivial product structure nearby. Thus, inside the plug  $M = L_0 \times \mD^q$, one can select a finite set of compact leaves without holonomy for $\cH_1$, and then apply the construction of this paper to these leaves individually.  Each of these modifications is restricted in its construction only by the fact that the  radii of the tubes used must decrease  by some proportion. The result is a   foliation $\cH_2$. Note that every leaf of $\cH_2$ is again a covering of the base leaf $L_0$. Now, $\cH_2$ is again a foliation which  contains open saturated sets, foliated by compact leaves with trivial product structure nearby. Thus, the process can be repeated for another collection of compact leaves without holonomy of $\cH_2$ to yield a foliation $\cH_3$. Repeat this process inductively, to obtain a limiting $C^0$-foliation $\cH_{\infty}$ which is distal and preserves a smooth transverse volume form. If the choices for the data used to construct each foliation $\cH_{\ell}$ are appropriately chosen, then estimates as used in section~\ref{sec-estimates}, show that $\cH_{\infty}$ will also be $C^r$ for a specified $1 \leq r \leq \infty$.
Note that if the number of leaves chosen at each stage is greater than one, then the resulting foliation $\cH_{\infty}$ will contain an uncountable number of solenoidal minimal sets. By the above remarks, one can even achieve that all of these solenoidal minimal sets with base $L_0$  are      pairwise non-homeomorphic.

\section{Remarks and Questions}\label{sec-remarks}

 In the following,    $L_0$ denotes a closed oriented connected manifold of
dimension $n \geq 1$, and $M = L_0 \times
\mD^q$ is the product disk bundle over $L_0$ for $q \geq 2$. Note that  $H^n(L_0 ; \mZ) \cong \mZ$.

 Markus and Meyer have shown \cite{MM1980} that the generic
Hamiltonian flow of a compact manifold of dimension at least four
contains as a limit set a one-dimensional solenoid from
\textit{each} homeomorphism class.  Can a similar result be true for
solenoids which occur as minimal sets of foliations?
\begin {quest}\label{quest-unctble}
Let $\whF$ be a $C^r$-foliation of a compact manifold $M$ with
leaves of dimension $n \geq 2$. Suppose that $L_0$ is a compact leaf
of $\F$ without holonomy. Does there exists a $C^r$-perturbation of
$\F$ such that $\F$ has an uncountable number of solenoidal minimal
sets, which realize every homeomorphism class of solenoid with base
$L_0$?
\end{quest}
If we restrict the question only to solenoids of product type over   $L_0$ which has abelian fundamental group, then the iterative construction above, and the result of  Markus and Meyer  \cite{MM1980} cited above,  suggests the answer may be yes. On the other hand, there are also uncountably many solenoids over $L_0$ which are not of product type, so the problem is also about whether   all of these     non-product  types can be realized by methods similar to those of this paper.

This leads directly to an important question, unresolved by the constructions using the standard representations in section~\ref{sec-foliations}.
\begin {quest}\label{quest-tori1}
Does every    solenoid $\cS$ with base $\mT^k$, for $k \geq 2$,
admit an embedding into a $C^{\infty}$-foliation $\whF$ of $M =
\mT^k \times \mD^q$, where $q \geq 2k$, so that its leaves are
covers of the zero section $L_0 \times \{0\} \subset M$? If not, is
there a natural invariant  of a solenoid  with base $\mT^k$, which
determines whether $\cS$ admits such an embedding?
\end{quest}
It seems plausible, based on our methods, that there may be
solenoids of dimension greater than one such that no presentation of $\cS$ can be embedded
smoothly. One approach to identify such an $\cS$ would be to find
some invariant of $\cS$ which implies that it is not possible to
find some presentation $\cP$ for which the inductive methods of this
paper yield a smooth embedding. Here is a related problem:

\begin {quest}\label{quest-tori2}
Let  $\cS$ be a solenoid with base $\mT^k$, for $k \geq 2$. Assume that $\cS$  admits an
 embedding into a $C^1$-foliation $\whF$ of $M$, where $q \geq 2k$, so that its leaves are covers of the zero section $L_0 \times \{0\} \subset M$. Is there some $1 \leq r \leq  \infty$ such that $\cS$ admits a $C^r$-embedding, but not a $C^{r+1}$-embedding?
\end{quest}
Finally, our methods leave the following question completely unresolved:
\begin {quest}\label{quest-tori3}
Let  $\cS$ be a solenoid with base $\mT^k$, for $k \geq 2$. Does
$\cS$ always admit an embedding in a codimension $q = 2$ foliation?
More generally,  for $q \geq 2$,  find an invariant of $\cS$ which guarantees that
$\cS$ admits an embedding into a codimension $q$ foliation.
\end{quest}
We show in this work,  that one can always embed $\cS$ into codimension $q = 2k_1$ where $k_1 \geq 1$ is the ``rank'' of the compact abelian Cantor group $\mK_0$ which is the fiber of $\cS \to \mT^k$.

Question ~\ref{quest-tori3} would be solved if we had a complete understanding of  the relation between the Cantor group fiber of a McCord solenoid and the codimension of possible embeddings. This suggests another line of questions, again equally unknown.

\begin {quest}\label{quest-McCord}
Let  $\cS$ be a McCord solenoid with base $L_0$ and Cantor group
fiber $\mK_0$. Does $\cS$ always admit an embedding into a $C^1$-foliation $\whF$ of $M$, so that its leaves are covers of the zero
section $L_0 \times \{0\} \subset M$? Are there invariants of the
algebraic or topological structure of the fiber $\mK_0$ which
determine whether such an embedding exists and the possible
codimension for such an embedding?
\end{quest}
For $C^2$-embeddings, we expect the answer to Question~\ref{quest-McCord} is negative:
\begin {conj}\label{conj-abelian}
Let  $\cS$ be a McCord solenoid with base $L_0$ and Cantor group fiber $\mK_0$.
If $\cS$   admits an embedding into a $C^2$-foliation $\whF$ of $M$, then $\mK_0$ admits a clopen neighborhood of the identity which is an abelian subgroup.
\end{conj}

\medskip

The above questions   primarily concern how the internal structure of a solenoid $\cS$ impact on the problem of finding an embedding  of $\cS$ as a minimal set of a $C^r$-foliation, for $1 \leq r \leq \infty$. The last set of questions concern the external geometry and topology of such an embedding.

In general there is a significant difference in the types of
possible minimal sets as the codimension varies, and there is a
change in how the minimal sets can occur. In~\cite{BM1995,
EThomas1973} it is shown that in any neighborhood of a
one-dimensional solenoid occurring as a minimal set of continuous
flow in $\mR^3$ there will be periodic orbit, while the same is not
true in $\mR^4$ as shown by Bell and
Meyer \cite{BM1995}. Thus, we do not expect to
generally find closed leaves in neighborhoods of higher dimensional solenoids embedded in foliation where the codimension is greater than $2$, although  it is an open problem to give constructions of such examples. 
In our plug, in any neighborhood of the minimal
solenoid there will be closed leaves. This leads to the following
question.
\begin {quest}\label{quest-smoothrestriction2}
If a solenoid $\cS$ embeds as a minimal set   of a codimension two
$C^r$-foliation, for $r \geq 1$, does every neighborhood of the embedded solenoid contain a compact leaf?

More generally, suppose the Cantor group fiber $\mK_0$ of a McCord solenoid $\cS$ has rank $k$,   is there  a   minimal codimension $q$, depending only on $k$, such that every embedding of $\cS$ has the ``compact leaf'' property?
\end{quest}
Note that a positive solution to Question~\ref{quest-smoothrestriction2}, or its modified version,  implies that
we can recover a presentation of $\cS$ from such an embedding. For $\ell \geq 1$, let   $U_{\ell}$ be the open neighborhood
of $\cS$ consisting of points of $M$ which are distance less than  $\e_{\ell} = 1/\ell$ from  $\cS$. Then choose a compact
leaf $L_{\ell} \subset U_1 \cap \cdots \cap U_{\ell}$, and the sequence $\{L_{\ell} \mid \ell \geq \ell_0\}$ is a presentation
for $\cS$ for $\ell_0$   sufficiently large.

\medskip
Our results show that any homeomorphism class of a one-dimensional
solenoid can be embedded as a minimal set of a $C^{\infty}$
foliation in $\mR^3$, but we have not completely answered
Problem~\ref{prob-embed2} even in the one-dimensional case, and so
it remains largely open. In~\cite{GST1994} restrictions are given
for how ``twisted" the presentation of a $C^1$ embedded presentation
of a one-dimensional solenoid can be in a solid torus. One can
supplement Problem~\ref{prob-embed2} by requiring various
cohomological conditions reflecting the twisting of the compact
leaves $L_i$, leading to an entire family of related problems.

%%%%%%%%%%%%%%%%%%%%%%%%%%%%%%%%

\end{document}